%% file: Hen2_main.tex
\documentstyle{amsppt}
\magnification\magstep1
\input imsmark
\input GlobalMacros
\nologo
\bigsinglespacing
\NoRunningHeads
\fulluntitledpage
\SetFiguresAnywhere % as opposed to \SetFiguresOnNewPage or \SetFiguresOnTop
\def\dotto{\underset\times\to\rightarrow}
%***************
\def \introchap {1}
   \def \IndLimFig{\introchap.1}
   \def \pcheckOK{\introchap.2}
   \def \diskstopoints{\introchap.3}
   \def \maintheorem{\introchap.4}
%******************
\def \telescopesection {2}
   \def \pCovers{\telescopesection.1}
   \def \teldefinespoint{\telescopesection.2}
   \def \telexample{\telescopesection.3}
   \def \strucstabpols{\telescopesection.4}
%*******************
\def \crosseschap {3}
   \def \crossfig{\crosseschap.1}
   \def \crossdef{\crosseschap.2}
   \def \degsofcrosses{\crosseschap.3}
   \def \manycurvesproper{\crosseschap.4}
   \def \horizdiskincross{\crosseschap.5}
   \def \intersectcurves{\crosseschap.6}
   \def \compsofcrosses{\crosseschap.7}
   \def \metricanddisks{\crosseschap.8}
   \def \conefieldsinv{\crosseschap.9}
   \def \conesnatural{\crosseschap.10}
   \def \onecrossescontract{\crosseschap.11}
   \def \unstablemansofcrosses{\crosseschap.12}
   \def \divergenceEnough{\crosseschap.13}
   \def \stabofcrosses{\crosseschap.14}
   \def \biinfinitesequences{\crosseschap.15}
%**********************************
\def \perturbpolschap {4}
   \def \crudepic{\perturbpolschap.1}
   \def \deltaneib{\perturbpolschap.2}
   \def \pictureofV{\perturbpolschap.3}
   \def \attcyclespersist{\perturbpolschap.4}
   \def \Visbidisk{\perturbpolschap.5}
   \def \howintersect{\perturbpolschap.6}
   \def \degoneremark{\perturbpolschap.7}
   \def \discwithsmalldiscs{\perturbpolschap.8}
   \def \structureofJ{\perturbpolschap.9}
   \def \whereJs{\perturbpolschap.10}
%********************************
\def \unstablemansection {5}
   \def \unstabmanthm{\unstablemansection.1}
   \def \improveremark{\unstablemansection.2}
   \def \projtoU{\unstablemansection.3}
   \def \existProjNearBound{\unstablemansection.4}
   \def \hyperplanes{\unstablemansection.5}
   \def \troublewithdiffeq{\unstablemansection.6}
%********************************
\def \stablemansection {6}
   \def \foliationofCCech{\stablemansection.1}
   \def \stabmantheorem{\stablemansection.2}
   \def \stabmantojulia{\stablemansection.3}
   \def \strucOfWs{\stablemansection.4}
   \def \homeoonpartescaping{\stablemansection.5}
   \def \infdimconst{\stablemansection.6}
   \def \homeopic{\stablemansection.7}
   \def \Hamstrom{\stablemansection.8}
   \def \quasiconformal{\stablemansection.9}
   \def \asymdev{\stablemansection.10}
   \def \easyextension{\stablemansection.11}
   \def \phipluscont{\stablemansection.12}
   \def \CChechWellDefined{\stablemansection.13}
%*******************************
\def \examplesection {7}
   \def \projnocrit{\examplesection.1}
   \def \BadAtAttCycle{\examplesection.2}
   \def \BadExRemark{\examplesection.3}
   \def \HomeoBetweenAnnuli{\examplesection.4}
   \def \JPlusisSphere{\examplesection.5}
   \def \AccessBdyEx{\examplesection.6}
   \def \boundaryOfBasin{\examplesection.7}
   \def \clamexample{\examplesection.8}
   \def \Hawaii{\examplesection.9}
   \def \homologyofbasin{\examplesection.10}
   \def \Jplussimplyconnected{\examplesection.11}
   \def \BigFundGroupRemark{\examplesection.12}
   \def \SingvCechRemark{\examplesection.13}
   \def \hornexample{\examplesection.14}
   \def \horn{\examplesection.15}
%*******************************
\def \Wadasection {8}
   \def \goodpols{\Wadasection.1}
   \def \enoughtofindpointsinT{\Wadasection.2}
   \def \periodninefig{\Wadasection.3}
   \def \renormalizabilityremark{\Wadasection.4}
   \def \wadalakestheorem{\Wadasection.5}
   \def \notcomplexgenerality{\Wadasection.6}
   \def \LakesA{\Wadasection.7}
   \def \LakesB{\Wadasection.8}
   \def \LakesC{\Wadasection.9}
   \def \Lakes{\Wadasection.10}
%*******************************
{
\input hen2_top.tex
\vfil\eject
}
\input hen2_intro.tex

\input hen2_telescopes.tex
\input hen2_crosses.tex
\input hen2_perturbpols.tex
\input hen2_unstabman.tex
\input hen2_stabman.tex
\input hen2_examples.tex
\input hen2_Wada.tex
\input hen2_ref.tex
\enddocument

%% file: imsmark.tex
\def\IMSmarkvadjust{0 pt}
\def\IMSmarkhadjust{0 pt}
\def\IMSmarkhpadding{0 pt}
\def\IMSpubltext{Published in modified form:}
\def\SBIMSMark#1#2#3{
 \font\SBF=cmss10 at 10 true pt
 \font\SBI=cmssi10 at 10 true pt
 \setbox0=\hbox{\SBF \hbox to \IMSmarkhpadding{\relax}
                Stony Brook IMS Preprint \##1}
 \setbox2=\hbox to \wd0{\hfil \SBI #2}
 \setbox4=\hbox to \wd0{\hfil \SBI #3}
 \setbox6=\hbox to \wd0{\hss
             \vbox{\hsize=\wd0 \parskip=0pt \baselineskip=10 true pt
                   \copy0 \break%
                   \copy2 \break% 
                   \copy4 \break}}
 \dimen0=\ht6   \advance\dimen0 by \vsize \advance\dimen0 by 8 true pt
                \advance\dimen0 by -\pagetotal
	        \advance\dimen0 by \IMSmarkvadjust
 \dimen2=\hsize \advance\dimen2 by .25 true in
	        \advance\dimen2 by \IMSmarkhadjust

%
%   Check for publication info
%
%  \newread\jref
  \openin2=publishd.tex
  \ifeof2\setbox0=\hbox to 0pt{}
  \else 
     \setbox0=\hbox to 3.1 true in{
                \vbox to \ht6{\hsize=3 true in \parskip=0pt  \noindent  
                {\SBI \IMSpubltext}\hfil\break
                \input publishd.tex 
                \vfill}}
  \fi
  \closein2
  \ht0=0pt \dp0=0pt
 \ht6=0pt \dp6=0pt
 \setbox8=\vbox to \dimen0{\vfill \hbox to \dimen2{\copy0 \hss \copy6}}
 \ht8=0pt \dp8=0pt \wd8=0pt
 \copy8
 \message{*** Stony Brook IMS Preprint #1, #2. #3 ***}
}

%% file: publishd.tex
{\it Real and Complex Dynamical Systems}, 
edit. B.~Branner and P.~Hjorth,
{\it NATO Adv. Sci. Inst. Ser.~C Math. Phys. Sci.}, {\bf 464},
Kluwer Acad. Publ., 1995.

%% file: hen2_top.tex
\topmatter
\title\nofrills\bigbf
H\'enon Mappings in the Complex Domain \\ II:  projective and inductive limits of polynomials
\endtitle
\author
John H. Hubbard and Ralph W. Oberste-Vorth
\endauthor
\address
Department of Mathematics, Cornell University, Ithaca, New York 14850
\endaddress
\email
hubbard \@ math.cornell.edu
\endemail
\address
Department of Mathematics, University of South Florida, Tampa, Florida 33620-5700
\endaddress
\email
ralph \@ math.usf.edu
\endemail
%\rightheadtext\nofrills{H\'enon Mappings in the Complex Domain II}
%\leftheadtext\nofrills{John H. Hubbard and Ralph W. Oberste-Vorth}
%\date
%Preprint: \today
%\enddate
\abstract
Let $H: \C^2 \to \C^2$ be the H\'enon mapping given by
$$
\bvec{x}{y} \mapsto \bvec{p(x) - ay}{x}.
$$
The key invariant subsets are $K_\pm$, the sets of points with bounded forward images, $J_\pm = \d
K_\pm$ their boundaries, $J = J_+ \cap J_-$, and $K = K_+ \cap K_-$. In this paper we identify the
topological structure of these sets when $p$ is hyperbolic and $|a|$ is sufficiently small, \ie, when
$H$ is a small perturbation of the polynomial $p$. The description involves projective and inductive
limits of objects defined in terms of $p$ alone. 
\endabstract

\SBIMSMark{1994/1}{February 1994}{}
\toc
\widestnumber\head{7}
\head1. Introduction                            \page{ 2}\endhead
\head2. Telescopes and Hyperbolic Polynomials   \page{ 7}\endhead
\head3. Crossed Mappings                        \page{ 9}\endhead
\head4. Perturbations of hyperbolic polynomials \page{15}\endhead
\head5. Characterization of $J_-$               \page{21}\endhead
\head6. Characterization of $J_+$               \page{25}\endhead
\head7. Examples                                \page{32}\endhead
\head8. Lakes of Wada in Dynamical Systems      \page{40}\endhead
\head {} References                             \page{46}\endhead
\endtoc
\endtopmatter

%% file: hen2_intro.tex
\heading
\introchap. Introduction
\endheading

This paper continues the study, begun with \cite{HO}, of the H\'enon family of mappings as a family of
mappings of two complex variables. Let $p(z)$ be a polynomial in one variable and $a \ne 0$ a complex
number. A H\'enon mapping is one which can be written
$$
H= H_{p,a}: \bmatrix x \\ y \endbmatrix \mapsto \bmatrix p(x) - ay \\ x \endbmatrix.
$$
Such a mapping has Jacobian $a$, and if $a \ne 0$, it is invertible:
$$
H^{-1}_{p,a}: \bmatrix x \\ y \endbmatrix \mapsto \bmatrix y \\ (p(y)-x)/a \endbmatrix.
$$

The key invariant subsets under such a mapping are
$$
K_\pm = \set {\bvec xy \in \C^2}{\left\|H^{\circ n} \bvec xy \right\| \text{ bounded as  } n \to \pm
\infty}
$$
as well as
$$
J_\pm = \d K_\pm, \quad K = K_+ \cap K_-, \quad \text{and} \quad J = J_+ \cap J_-.
$$

When $a = 0$, the degenerate H\'enon mapping
$$
H_{p,0}: \bmatrix x\\y \endbmatrix \mapsto \bmatrix p(x) \\ x \endbmatrix
$$
is not invertible, but maps all of $\C^2$ to the curve $C_p$ of equation $x = p(y)$, and reduces to $x
\mapsto p(x)$ in the first coordinate.

According to the theory that hyperbolic dynamics is stable under perturbations, you would expect that
$H_{p,a}$ could be understood as a perturbation of $p$ for $a$ sufficiently small when $p$ is
hyperbolic. Many people (\eg, Holmes, Whitley, and Williams, \cf, \cite{Ho}, \cite{HWh}, and
\cite{HWi} for further references) have done this in the real domain, at least for $|a|$ small.
Benedicks and Carleson have gone further in this direction \cite{BC}. In this article we will do the
same in the complex domain. By the techniques used here we can only deal with perturbations of
hyperbolic polynomials, and not the much more difficult ones studied by Carleson and Benedicks.

There is a fundamental conflict between the H\'enon mapping and polynomials: polynomials are not
injective and H\'enon mappings are. We will describe two ways of creating from a polynomial $p$
objects which do carry bijective dynamics; both appear as invariant subsets of $\C^2$ for H\'enon
mappings which are sufficiently small perturbations of hyperbolic polynomials.

\subheading{The projective limit construction}

Let
$$
\hat \C_p = \varprojlim (\C, p).
$$
A point of this projective limit is a point $z_0 \in \C$ and a {\it history of the point} $z_0$ under
the iteration of $p$. More precisely,
$$
\hat{\C}_p = \bigl\{\, (\dots, z_{-2}, z_{-1}, z_0) \bigm| p(z_{-i-1}) = z_{-i} \text{ for all } i =
\dots, -2,-1,0 \bigr\}.
$$
The mapping $p$ induces a mapping $\hat p: \hat{\C}_p \to \hat{\C}_p$ by
$$
\hat p(\dots, z_{-2}, z_{-1}, z_0 ) = (\dots, p(z_{-2}), p(z_{-1}), p(z_0)) = (\dots, z_{-1}, z_0,
p(z_0))
$$
which is of course bijective:
$$
\hat p^{-1}(\dots, z_{-2}, z_{-1}, z_0) = (\dots, z_{-2}, z_{-1} ).
$$
In section \examplesection, we will give a description of this space which makes it reasonably
understandable when $p$ is hyperbolic.

\subheading {The inductive limit construction}

Recall that if $f: X \to X$ is a mapping from a space to itself, then the inductive limit
$$
\check X_f = \varinjlim (X,f)
$$
is the quotient $(X \times \N)/\sim$, where $\sim$ is generated by setting $(x,n) \sim (f(x),n+1)$.

%\setEPSF EPSF.IndLimFig 4.69in by 1.64in caption (Figure \IndLimFig:
%Inductive limit as an increasing union)
\setPSFig 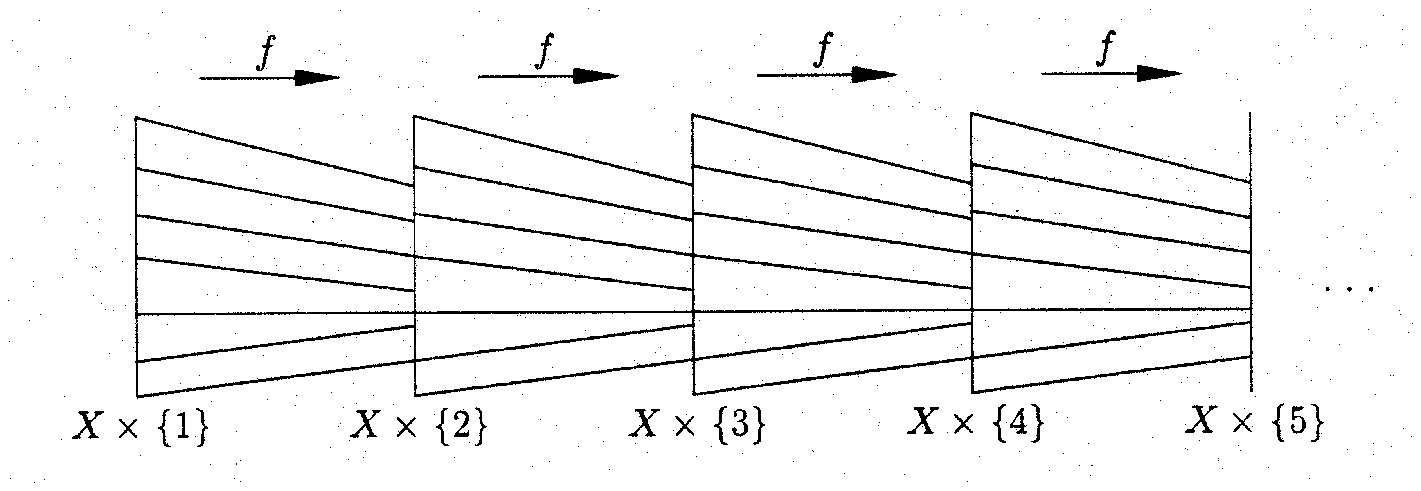 {2true in} caption (Figure \IndLimFig:
 Inductive limit as an increasing union)

Inductive limits are pathological objects in general, and will be Hausdorff only when $f$ has some
nice properties. We will consider them only when $f$ is open and injective, in which case the
inductive limit is an increasing union of subsets homeomorphic to $X$, hence locally as nice as $X$.

The inductive limit comes with a map to itself: $\check f: \check X_f \to \check X_f$ induced by
$$
\check f: (x,n) \mapsto (f(x),n) \sim (x, n-1).
$$
This mapping is obviously bijective, as an inverse is induced by $(x,n) \mapsto (x,n+1)$.

We will now apply this construction to polynomials.  Our construction only makes sense for
polynomials $p$ with no critical points in the Julia set; however, we will only apply it to hyperbolic
polynomials, which all have this property. Let $D \subset \C$ be a disk of radius $R$ sufficiently
large so that $J_p \subset D$, where $J_p$ is the Julia set of $p$. 

Consider the mapping $f_{p,\alpha,R}: J_p \times D \to J_p \times \C$ given by 
$$
f_{p,\alpha,R}(\zeta,z) = \left(p(\zeta), \zeta + \alpha \frac z{p'(\zeta)}\right), 
$$
which is well defined since $p'(\zeta) \ne 0$.  

\proclaim {Lemma \pcheckOK}
If $p$ is hyperbolic, and if $|\alpha| \ne 0$ is sufficiently small, then the image $f_{p,\alpha,R}
(J_p \times D)$ is contained in $J_p \times D$ and $f_{p,\alpha,R}$ is open and injective. 
\endproclaim

\demo {Proof}
Recall that if $p$ is hyperbolic, there are no critical points of $p$ in $J_p$ (in fact, this is the
only property of hyperbolic polynomials this lemma requires). Thus the formula is well defined, and
clearly if $|\alpha|$ is sufficiently small, the image lies in $J_p \times D$. Moreover, if there are
no critical points in $J_p$, then there exists $\epsilon >0$ such that when $\zeta_1 \ne \zeta_2 \in
J_p$ and $p(\zeta_1) = p(\zeta_2)$, then $|\zeta_1 - \zeta_2| > \epsilon$. If we choose
$$
0 < |\alpha | < \frac {\epsilon R}{\inf\limits_{\zeta\in J_p}\, |p'(\zeta)|},
$$
then $f_{p,\alpha,R}$ is clearly injective. The mapping is open because it is a local homeomorphism.
\QED
\enddemo

Thus when $p$ is hyperbolic and $|\alpha|$ is sufficiently small and $R$ is sufficiently large, we may
set
$$
\check \C_p = \check \C_{p,\alpha,R} = \varinjlim (J_p \times D, f_{p,\alpha,R}),
$$
and we will denote by 
$$
\check p = \check f_{p,\alpha,R}:\check \C_p \to \check \C_p
$$
the bijective mapping above.

If $\psi:X \to Y$ is a homeomorphism conjugating $f:X \to X$ and $g:Y \to Y$, then $\psi$ induces a
homeomorphism $\check \psi:\check X_f \to \check Y_g$ conjugating $\check f: \check X_f \to \check X_f$
to $\check g: \check Y_g \to \check Y_g$. Thus the following proposition, which is proved in Section
\stablemansection, shows that we can drop the indices $\alpha$ and $R$, and speak simply of $\check
p: \check \C_p \to \check \C_p$.

\proclaim {Proposition \CChechWellDefined}
For all $\alpha_1,\alpha_2$ sufficiently small and all $R_1$ and $R_2$ sufficiently large, there is a
homeomorphism 
$$
\psi: J_p \times D_{R_1} \to J_p \times D_{R_2}
$$
conjugating $f_{p,\alpha_1,R_1}$ to $f_{p,\alpha_2,R_2}$.
\endproclaim 

This justifies writing simply $f_p$ and $\check \C_p$. The space $\check \C_p$ is quite difficult to
understand. The only case where it is anything familiar is when $J_p$ is a Jordan curve; in that case
$\check \C_p$ is homeomorphic to the complement of a solenoid in a $3$-sphere. Proposition
\foliationofCCech\ gives some important information, and much more is shown in Section
\examplesection. In Section \Wadasection, we show that when $p$ is a real hyperbolic polynomial, the
real part $\check \R_p$ is often the common separator of Lakes of Wada. This illustrates some of the
unavoidable complexity.

\subheading {An embedding of $\hat J_p$ into $\check \C_p$}

The inductive and projective limits above are related: the projective limit $\hat J_p$ is naturally an
invariant subset of both. This is obvious for $\hat \C_p$; let us see why it is true for $\check
\C_p$.

Let $\underline \zeta = (\dots, \zeta_{-2}, \zeta_{-1}, \zeta_0) \in \hat J_p$, and consider the
intersection
$$
(\{\zeta_0\} \times D) \,\cap\, \check p(\{\zeta_{-1}\}\times D) \,\cap\, \cdots \,\cap\, \check
p^{\circ n} (\{\zeta_{-n}\} \times D) \,\cap\, \cdots
$$
\proclaim {Lemma \diskstopoints}
This is a nested sequence of embedded disks, and the intersection is a single point.
\endproclaim

\demo {Proof}
The nesting is obvious. As we have defined it, there exists a disk $D_1$ relatively compact in $D$
such that
$$
f_{p,\alpha,R}(J_p \times D) \subset J_p \times D_1.
$$
There are infinitely many disjoint conformal copies of the annulus $D \setminus \bar D_1$ surrounding
the intersection above. This shows that the intersection is a point.
\QED
\enddemo

Let us call $\psi: \hat J_p \to \check \C_p$ the mapping which associates to $\underline \zeta$ the
unique point in the above intersection. Clearly the diagram
$$
\CD
  \hat J_p   @>\psi>>   \check \C_p    \\
@V\hat p VV            @VV \check p V  \\
  \hat J_p   @>>\psi>   \check \C_p
\endCD
$$
commutes.

We will see in section \examplesection\ some examples of the objects above. In particular, we will see
that the construction above corresponds to seeing the solenoid as a projective limit of circles or a
decreasing intersection of solid tori.

\subheading {Riemann surface laminations}

It is rather difficult to find any category to which $\hat \C_p$ and $\check \C_p$ belong. A
first attempt is to say that they are are (or have large subsets which are) {\it Riemann
surface laminations}. For future reference, we define this category to have: 
\newline 
\indent{\it Objects:}
Hausdorff spaces which are locally products of Riemann surfaces by topological spaces, glued
together by local isomorphisms;  
\newline
\indent{\it Morphisms:}
Continuous mappings, analytic on each Riemann surface.

You should imagine the topological factor to be like a Julia set, either a Cantor Julia set (for
$\hat \C_p$) or a connected Julia set (for $\check \C_p$). This category has recently turned up in
several fields, and Sullivan's paper \cite{S} contains some basic material about this category.
Pictures of the H\'enon attractor \cite{H\'e2} or of basin boundaries will show that such
structures should be relevant to dynamical systems.

\subheading {The main result}

Both of the constructions above give objects which arise in the dynamical plane $\C^2$ of H\'enon
mappings.

\proclaim {Theorem \maintheorem}
Let $p$ be a hyperbolic polynomial. There exists $A$ such that if $0 < |a| < A$, then there exist
homeomorphisms
$$
\Phi_-: \hat \C_p \to J_- \qquad \text{ and } \qquad \Phi_+: \check \C_p \to J_+
$$
such that the diagrams
$$
\CD
 \hat \C_p   @>\Phi_->>    J_-    \\
@V\hat p VV              @VV H V  \\
 \hat \C_p   @>>\Phi_->    J_-
\endCD
\qquad \text{ and } \qquad
\CD
 \check \C_p   @>\Phi_+>>    J_+    \\
@V\check p VV              @VV H V  \\
 \check \C_p   @>>\Phi_+>    J_+
\endCD
$$
commute. On $\hat J_p$, the mappings $\Phi_+$ and $\Phi_-$ coincide, \ie, we have
$$
\Phi_- \mid_{\hat J_p} = \Phi_+ \circ \psi.
$$
\endproclaim

\subheading {Outline of the paper}

The proof we will give of Theorem \maintheorem\ is an adaptation of the technique of {\it
telescopes}, which we learned from Sullivan many years ago.

In section \telescopesection, we will review Sullivan's construction. This will serve several
purposes: it will motivate our construction, it will provide us with some constructions which we need,
and it will provide a written account of Sullivan's proof, which was never published.

In section \crosseschap, we will define our 2-dimensional analogs of expanding maps, which we call
{\it crossed mappings}. It seems clear that these are going to be of interest in many other settings,
and we have proved the basic results concerning them with considerable care.

In section \perturbpolschap, we show that for H\'enon mappings which are small perturbations of
hyperbolic polynomials, the mappings analogous to the telescope mappings are crossed mappings. This
will give us a homeomorphism $\Phi: \hat J_p \to J$ conjugating the H\'enon mapping to $\hat p$, and
locally the stable and unstable manifolds will also drop out of the construction.

In section \unstablemansection, we identify the unstable manifold of $J$ with $\hat \C_p$, and in
section \stablemansection\ we identify the stable manifold with $\check \C_p$. This last step is quite
delicate, and is surely the hardest proof in the paper.

Finally, in sections \examplesection\ and \Wadasection, we show in some examples exactly what these
results give us for the topology of H\'enon mappings, including Lakes of Wada.

%% file: hen2_telescopes.tex
\heading
\telescopesection. Telescopes and Hyperbolic Polynomials
\endheading

Many years ago, we learned from Sullivan that hyperbolic polynomials (and rational functions) are
structurally stable on their Julia sets. Sullivan used {\it telescopes} in his proof, and we are
planning to adapt this construction to H\'enon mappings.

Let $p(z)$ be a hyperbolic polynomial. In fact, everything we will say goes over to rational
functions without modification. We will take as our definition of hyperbolic that all critical
points are attracted to attractive periodic cycles. As we will see below, this is equivalent to
saying that $p$ is strongly expanding on the Julia set $J_p$.

Call $\Omega$ the Fatou set of $p$, $C$ the set of attracting periodic points of $p$ (including
$\infty$ ), and
$$
X_0 = \set {z \in \Omega}{d_\Omega(z,C)\le 1},
$$
where $d_\Omega$ is the Poincar\'e metric on $\Omega$. (The number 1 in the definition is arbitrary;
everything would go through for any positive constant.)

Further set $X_n = p^{-n}(X_0)$. The $X_n$ form an increasing collection of compact subsets of
$\Omega$ which exhaust $\Omega$, and they are strictly increasing in the sense that $X_{n-1}$ is
contained in the interior of $X_n$. Similarly, the sets $U_n = \overline \C \setminus X_n$ form a
basis of nested open neighborhoods of $J_p$, each relatively compact in the previous.

Let $N$ be the smallest index such that all the critical points of $p$ are in $X_N$.  Such an $N$
exists since there are only finitely many critical points, all in $\Omega$.

\proclaim {Proposition \pCovers}
The mapping $p: U_n \to U_{n-1}$ is a covering map for $n \ge N$. In particular, it is strongly
infinitesimally expanding for the Poincar\'e metric of $U_{n-1}$.
\endproclaim

\demo {Proof}
Clearly $p: U_n \to U_{n-1}$ is proper and a local homeomorphism, hence a covering map and a local
isometry from the Poincar\'e metric of $U_n$ to the Poincar\'e metric of $U_{n-1}$. Since $U_n$ is
relatively compact in $U_{n-1}$, the inclusion is strongly contracting for the Poincar\'e metric of
$U_{n-1}$.
\QED
\enddemo

We will call $U = U_{n-1}$ and $U'= U_n$, so that $p: U'\to U$ is a covering map. Choose $\epsilon >
0$ sufficiently small that for any $z \in U$, the set
$$
U_z = \set {z_1 \in U}{d_U(z_1,z) <\epsilon}
$$
is homeomorphic to a disk, and that $p$ restricted to $U_z$ is a homeomorphism to its image.

For any $z \in J_p$, define $U^0_z = U_z$, and recursively set
$$
U_z^n = U_z^{n-1} \,\cap\, p^{-n}(U_{p^{\circ n}(z)}).
$$
It is easy to show that each $U_z^n$ is homeomorphic to a disk.

\proclaim {Proposition \teldefinespoint}
We have
$$
\{z\} = \bigcap_n\, U_z^n.
$$
\endproclaim

\demo{Proof}
Clearly $z$ is in the intersection; the only problem is to show that the intersection is a single
point. This follows from the strong expansion: if $p$ expands by a factor of $K >1$, then the
diameter of $U_z^n$ is at most $\epsilon/K^n$.
\QED
\enddemo

Sullivan defines a {\it $p$-telescope}\/ to be a sequence of disks $W_0, W_1, \dots$ such that
$W_{n+1}$ is relatively compact in $p(W_n)$.

\demo{Example \telexample}
If $z \in J_p$, the sequence of disks $U_z, U_{p(z)}, \dots$ is a telescope, and Proposition
\teldefinespoint\ says that a telescope defines a point. But clearly a telescope for $p$ is also a
telescope for a small perturbation of $p$, so that going from points to telescopes to points provides
a conjugacy between the Julia set of a hyperbolic polynomial and that of a small perturbation. This
is the idea behind Sullivan's proof.
\enddemo

\proclaim {Theorem \strucstabpols}
For any neighborhood $V$ of the Julia set of $p$, there exists a neighborhood of $p$ in the
$C^1$-topology such that any $p_1$ in that neighborhood is conjugate to $p$ on a neighborhood of the
Julia set.
\endproclaim

\demo {Sketch of proof}
Define $\phi: J_p \to V$ by
$$
\phi(z) = \bigcap_n\, p_1^{-n}(U_{p^{\circ n}(z)}).
$$
Just as above, for $p_1$ sufficiently close to $p$ in the $C^1$-topology, this intersection is a
single point. A similar construction gives an inverse for $\phi$ on
$$
J_{p_1} = \set {z \in V}{p_1^{\circ n}(z) \in V \text{ for all $n$}}.
$$
Thus $\phi: J_p \to J_{p_1}$ is a homeomorphism conjugating $p$ to $p_1$ on the Julia sets. We
leave to the reader to verify that this homeomorphism can be extended to a neighborhood of $J_p$
which still conjugates $p$ to $p_1$.
\QED
\enddemo

%% file: hen2_crosses.tex
\heading
\crosseschap. Crossed Mappings
\endheading

In one dimension, the mappings useful for structural stability are those which map a disk strictly
outside another. In higher dimensions, we will be interested in bijective mappings defined on
bidisks which map the ``horizontal boundary'' outside itself, and the inverses of which map
the ``vertical boundary'' outside itself. Thus, they ``look like'' Figure \crossfig.

%\setEPSF EPSF.crossfig 2.97in by 1.03in caption (Figure \crossfig: A
%$1$-crossed mapping)
\setPSFig 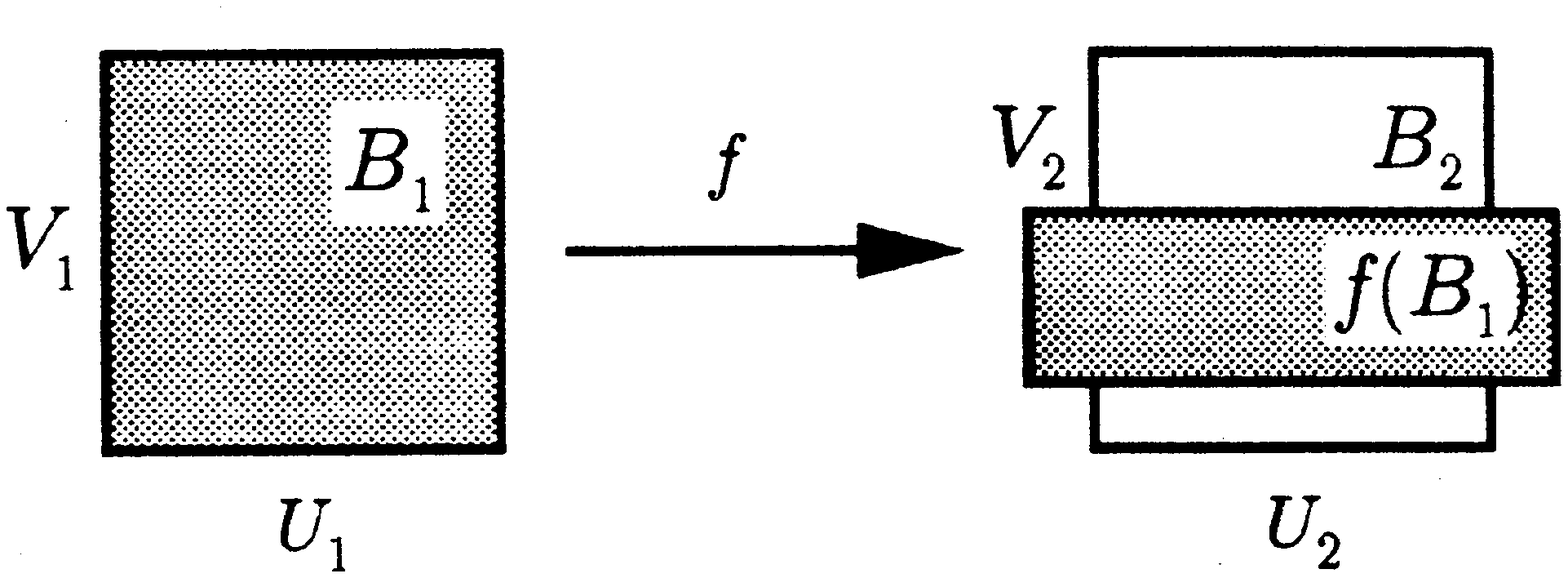 {3.5 true cm} caption (Figure \crossfig: A $1$-crossed mapping)

We need to formalize what this means. Let $B_1 = U_1 \times V_1$ and $B_2 = U_2 \times
V_2$ be bidisks.

\definition {Definition \crossdef}
A crossed mapping from $B_1$ to $B_2$ is a triple $(W_1, W_2, f)$, where
\roster
\item
$W_1 \subset U'_1 \times V_1$ where $U'_1 \subset U_1$ is a relatively compact open subset,
\item
$W_2 \subset U_2 \times V'_2$ where $V'_2 \subset V_2$ is a relatively compact open subset,
\item
$f: W_1 \to W_2$ is a holomorphic isomorphism, such that for all $y \in V_1$, the mapping
$$
\pr_1 \circ f\vert_{W_1\cap (U_1 \times \{y\})}: W_1 \cap (U_1 \times \{y\}) \to U_2
$$
is proper, and the mapping
$$
\pr_2 \circ f^{-1}\vert_{W_2\cap(\{x\}\times V_2)}: W_2 \cap (\{x\}\times V_2) \to V_1
$$
is proper.
\endroster
\enddefinition

To make the notation less cumbersome, we will often write $f: B_1 \dotto B_2$ for a crossed mapping,
leaving the precise $W_1$ and $W_2$ to be determined by the context.

\proclaim {Proposition \degsofcrosses}
If $f: W_1 \to W_2$ is a crossed mapping from $B_1$ to $B_2$, then all maps
$$
\gather
\pr_1 \circ f\vert_{W_1\cap (U_1 \times \{y\})}: W_1 \cap (U_1 \times \{y\}) \to U_2 \\
\intertext{and}
\pr_2 \circ f^{-1}\vert_{W_2\cap(\{x\}\times V_2)}: W_2 \cap (\{x\}\times V_2) \to V_1
\endgather
$$
have the same degree, which will be called the {\it degree}\/ of the crossed mapping.
\endproclaim

\demo {Proof}
Choose $x \in U_2$, and consider $Z_x = f^{-1}\left(W_2 \cap (\{x\} \times V_2)\right)$, which is a
closed analytic curve in $B_1$ (\ie, a Riemann surface closed in $B_1$). The mapping $\pr_2: Z_x \to
V_1$ is proper, hence a finite ramified covering map, of some degree $k(x)$. For every $y \in V_1$,
the line $U_1 \times \{y\}$ cuts $Z_x$ in precisely $k(x)$ points, counted with multiplicity (where
the multiplicity is almost by definition the local degree of the projection above). But for each such
$y$, these $k(x)$ intersection points are mapped by $f$ exactly to the intersections of $f(W_1 \cap
(U_1 \times \{y\}))$ with the line $\{x\} \times V_2$; these count the degree of
$$
\pr_1 \circ f\vert_{W_1 \cap (U_1 \times \{y\})}: W_1 \cap (U_1 \times \{y\}) \to U_2.
$$
Thus these maps all have the same degree, and the same argument applied to $f^{-1}$ shows that the
maps
$$
\pr_2 \circ f^{-1}\vert_{W_2 \cap (\{x\} \times V_2)}: W_2 \cap (\{x\} \times V_2) \to V_1
$$
also all have the same degree; \ie, $k(x)$ does not depend on $x$. It is clear from the proof that
the two classes of mappings have the same degree.
\QED
\enddemo

Figure \crossfig\ represents a $1$-crossed mapping.

\proclaim {Proposition \manycurvesproper}
If $f: W_1 \to W_2$ is a crossed mapping from $B_1$ to $B_2$ and $X \subset B_1$ is an analytic curve
such that $\pr_1: X \to U_1$ is proper of degree $l$, then $\pr_1 \circ f: X \cap W_1 \to U_2$ is
proper of degree $kl$.
\endproclaim

\remark {Remark \horizdiskincross}
The case where $X$ is a horizontal disk is part of the definition of a crossed mapping.
\endremark

We require the following lemma, which is classical.

\proclaim {Lemma \intersectcurves}
Let $X$ and $Y$ be curves in a bidisk $B = U \times V$ such that $\pr_1: X \to U$ and $\pr_2: Y \to
V$ are proper of degrees $k_X$ and $k_Y$, respectively. Moreover, suppose that $\pr_1(Y)$ is
relatively compact in $U$. Then $X$ and $Y$ intersect in $k_Xk_Y$ points counted with multiplicity.
\endproclaim

\demo {Proof of Proposition \manycurvesproper}
For each $x \in U_2$, the curve $X$ and the curve 
$$
Y_x = f^{-1}\left(W_2 \cap (\{x\} \times V_2)\right)
$$
 satisfy
the hypotheses of Lemma \intersectcurves. So these curves intersect in $kl$ points independent of
$x$. But this means that every vertical line in $B_2$ intersects $f(X \cap W_1)$ in the same number
of points. Since $f(X \cap W_1)$ is clearly closed in $W_2$, this shows that it maps by a proper map
to $U_2$.
\QEDL{(Proposition~\manycurvesproper)}
\enddemo

\proclaim {Proposition \compsofcrosses}
(a) Let $f: W_1 \to W_2$ be a crossed mapping from $B_1$ to $B_2$ of degree $k$. Then $f^{-1}: W_2
\to W_1$ is also a crossed mapping if all the coordinates are flipped. \newline
(b) If $B_1$, $B_2$, and $B_3$ are bidisks, $W_1 \subset B_1$, $W_2 \subset B_2$,
$\tilde W_2 \subset B_2$, and $\tilde W_3 \subset B_3$, and $f_1: W_1 \to \tilde W_2$ and $ f_2:
 W_2 \to \tilde W_3$ are $k_1$- and $ k_2$-crossed mappings, then
$$
 f_2 \circ f_1: W_1 \cap f_1^{-1}(W_2) \to \tilde W_3 \cap  f_2(\tilde W_2)
$$
is a $k_1 k_2$-crossed mapping from $B_1$ to $B_3$.
\endproclaim

\demo{Proof}
Part (a) is obvious. For (b), observe that the sets $S_1 = W_1 \cap f_1^{-1}(W_2)$ and $S_2 =
\tilde W_3 \cap f_2(\tilde W_2)$ clearly satisfy conditions (1) and (2) of the definition; it remains
to show (3). For any $y \in V_1$, the curve $X_y = f_1\left(W_1 \cap (U_1 \times \{y\})\right)$
satisfies the hypothesis of Proposition \manycurvesproper, with respect to the crossed mapping
$f_2:  W_2 \to \tilde W_3$. So the projection $\pr_1: f_2(W_2 \cap X_y) \to U_3$ is
proper. Proposition \manycurvesproper\ also shows that this proper projection has degree $k_1 k_2$.
\QED
\enddemo

A bidisk $B = U \times V$ carries, like all bounded domains, the Kobayashi metric, which in this
case is easy to describe: it is the product of the Poincar\'e metrics of $U$ and of $V$.
Crossed mappings of degree $1$ have special expansion and contraction properties with respect to
this metric. If $\xi \in T_xU$, we will denote by $|(x,\xi)|_U$ the length of the
tangent vector for the (infinitesimal) Poincar\'e metric of $U$. 

A first observation about this metric is the following:

\proclaim {Lemma \metricanddisks}
A tangent vector $(\xi, \eta) \in T_{(x,y)}B$ is tangent to a disk which is the graph of an injective
mapping $g: U \to V$ if and only if $|(x,\xi)|_U \ge |(y,\eta)|_V$. This mapping can be taken to have
relatively compact image in $V$ if and only if $|(x, \xi)|_U > |(y,\eta)|_V$.
\endproclaim

\demo{Proof}
In one direction, this is Schwarz's Lemma: such a $g$ contracts in the Poincar\'e metrics and
contracts strictly if the image is relatively compact. In the other direction, by a biholomorphic
isomorphism, we may suppose the bidisk is the standard bidisk $D \times D$, and that $(x,y) =
(0,0)$. Then the line containing $(\xi, \eta)$ intersects the bidisk in an appropriate graph.
\QED
\enddemo

This gives us the appropriate tool to study crossed mappings of degree 1. For any bidisk $U
\times V$, consider the horizontal cone field $C_{(x,y)} \subset T_{(x,y)}B$ defined by
$$
C_{(x,y)} = \altset {(\xi,\eta) \in T_{(x,y)}B}{|(x,\xi)|_U \ge |(y,\eta)|_V}.
$$
Reversing the inequality gives the vertical cone field.

\proclaim{Proposition \conefieldsinv}
Let $f: W_1 \to W_2$ be a crossed mapping from $B_1$ to $B_2$ of degree 1. Then for all $(x,y) \in
W_1$, we have
$$
d_{(x,y)}f(C_{(x,y)}) \subset C_{f(x,y)}.
$$
Moreover, if $(\xi_1,\eta_1) \in C_{(x_1,y_1)}$, $f(x_1,y_1) = (x_2,y_2)$ and
$d_{(x_1,y_1)}f(\xi_1,\eta_1) = (\xi_2,\eta_2)$, then $|(x_2,\xi_2)|_{U_2} > |(x_1,\xi_1)|_{U_1}$.
\endproclaim

\remark{Remark \conesnatural}
This proposition illustrates the principle that in complex analysis,
inequalities often follow from topology.  In the complex setting, the existence of an invariant
cone-field is automatic; in the real it needs to be verified in each case \cite {Yoc1}, this is often
quite difficult.
\endremark 

\demo{Proof of Proposition \conefieldsinv}
A vector in the cone $C_{(x_1,y_1)}$ is tangent to a curve $X$ in $B_1$ proper of degree 1 over
$U_1$. The curve $f(X \cap W_1)$ is then proper of degree 1 over $U_2$ by Proposition
\manycurvesproper. Thus the tangent to $f(X \cap W_1)$ at $f(x_1,y_1)$ is in the cone
$C_{(x_2,y_2)}$. This proves the first part.

For the second part, observe that $\pr_1\circ f: X \cap W_1 \to U_2$ is an isomorphism, so that
$$
\pr_1 \circ \left(\pr_1\circ f\vert_{X \cap W_1}\right)^{-1}: U_2 \to U_1
$$
is an analytic mapping with relatively compact image. Thus it strictly contracts Poincar\'e lengths,
and its derivative maps $\xi_2$ to $\xi_1$.
\QEDL{(Proposition~\conefieldsinv)}
\enddemo

We will refer to smooth curves and surface in a bidisk $B=U \times V$ as {\it horizontal-like} of
{\it vertical-like} if their tangent spaces are in the horizontal or vertical cone respectively at
each of their points.

Suppose $B_1$, \dots, $B_{n+1}$ are bidisks such that $B_i = U_i\times V_i$. Suppose also that
$W_i\subset B_i$ ($i = 1, \dots, n$) and $\tilde W_i\subset B_i$ ($i=2, \dots, n+1$) are open
subsets so that $f_i: W_i \to \tilde W_{i+1}$ are crossed mappings of degree 1. Let
$$
S_1 = W_1 \cap f_1^{-1}(W_2) \dots \cap (f_1^{-1} \circ \dots \circ f_{n-1}^{-1})(W_{n-1})
$$
and
$$
S_2 = \tilde W_{n+1} \cap f_n(\tilde W_n) \cap \dots \cap (f_n \circ \dots \circ f_2)(\tilde W_2)
$$
so that Proposition \compsofcrosses\ and an obvious induction shows that the restriction $g$ of
$f_n \circ \dots \circ f_1$ to $S_1$ makes $g: S_1 \to S_2$ a crossed mapping of degree 1 from
$B_1$ to $B_{n+1}$.

We want to say that horizontal disks in $B_1$ intersect $S_1$ in regions with small diameter. It is
a bit harder to do this from Proposition \conefieldsinv\ than one might expect: the contraction
there is infinitesimal, and a domain $U' \subset U$ may have small diameter in $U$ although
the shortest curve joining points of $U'$ in $U'$ may still be long. We will take a different
tack, using complex analysis and moduli of annuli. Note that such methods, in a more
complicated context, have had great success recently in one-dimensional dynamics \cite{BH}, \cite{H},
\cite{Yoc2}.

Let $U$ be a simply connected Riemann surface isomorphic to the disk, and $U'$ a relatively
compact open subset. Define the {\it size of $U'$ in $U$}\/ to be the $1/M$, where $M$ is the
largest modulus of an annulus separating $U'$ from the boundary of $U$. We will really be interested
in the case where $U'$ is connected and simply connected; so that the size of $U'$ in $U$ is the
inverse of the modulus of $U-\overline{U'}$.

Note that the size is related to the Poincar\'e diameter by a double inequality; after conformal
mapping of $U$ to the disk, the extreme cases for a given diameter are a line segment and a round
disk.

\proclaim {Proposition \onecrossescontract}
Suppose that the size of the projection $\pr_1(W_i)$ in $U_i$ is $1/M_i$. Then for any $y \in V_1$,
the size in $U_1$ of $S_1 \cap (U_1 \times \{y\})$ is at least
$$
1 \left/ \ \sum_1^{n-1} \frac 1{M_i}\right..
$$
\endproclaim

\demo {Proof}
Consider the subsets
$$
W_1^j =f_1^{-1}\circ\dots \circ f_j^{-1}(W_{j+1}),
$$
which are nested so that
$$
B_1 \supset W_1 = W_1^1 \supset W_1^2 \supset \dots \supset W_1^n = S_1.
$$
For any $y \in V_1$, the annuli $(V_1\times\{y\})\cap (W_1^j \setminus W_1^{j+1})$ are disjoint
nested annuli for $j=1,\dots,n-1$, and the $j$th maps by $\pr_1\circ f_j\circ\dots\circ f_1$ to an
annulus which contains $U_{j+1} \setminus U'_{j+1}$, hence has modulus at least $1/M_{j+1}$. The
result now follows from the additivity of moduli.
\QED
\enddemo

\proclaim {Corollary \unstablemansofcrosses}
Let $B_0 = U_0\times V_0, B_1 = U_1\times V_1, \dots$ be an infinite sequence of bidisks, and
$f_i: B_i \dotto B_{i+1}$ be crossed mappings of degree 1, with  $U'_i$ of uniformly bounded size
in $U_i$. Then the set
$$
W^S_{\{B_n,f_n\}}= \set {\bmatrix x\\y\endbmatrix \in B_0}{ f_n \circ \dots \circ f_0\left(\bmatrix
x\\y\endbmatrix\right) \in B_n \quad \text {for all $n$}}
$$
is a vertical-like analytic disk in $B_0$, which maps by $\pr_2$ isomorphically to $V_0$, which we
will call the stable disk of the sequence of crossed mappings.
\endproclaim

Similarly, when we have backwards sequence of crossed mappings
$$
\dots \dotto B_{-1} \dotto B_0
$$
with uniformly bounded sizes, it will have a {\it unstable disk}, which will be horizontal-like.

\remark {Remark \divergenceEnough}
Rather than requiring that the sizes $1/M_i$ of the $U'_i$ in $U_i$ be uniformly bounded, it would be
enough to require that $\sum M_i = \infty$.
\endremark 

\demo{Proof of Corollary \unstablemansofcrosses}
For any $u_m \in U_m$, we can consider the set
$$
\Gamma_m =\set {(x,y) \in B_0}{f_{m-1} \circ \dots \circ f_0(x,y) \in \{u_m\} \times V_m}.
$$
This is a vertical-like analytic disk, so that there exists an inverse $\gamma_m: V_0 \to B_0$ of
$\pr_2$ which parameterizes it.

If $u_0, u_1, \dots$ is any sequence with $u_m \in U_m$, and $\gamma_m: V_0 \to B_0$ is constructed
as above for each $m$, then Proposition \onecrossescontract\ says that the sequence $\gamma_m$
converges uniformly. Clearly the limit is a parameterization of an vertical-like analytic disk
contained in $W^S_{\{B_n,f_n\}}$. Clearly by Corollary \unstablemansofcrosses, it is all of
$W^S_{\{B_n,f_n\}}$.
\QEDL{(Corollary~\unstablemansofcrosses)}
\enddemo

We will refer to $W^S_{\{B_n,f_n\}}$ as the stable set of the sequence of crossed mappings
$$
B_0 \overset f_0 \to {\dotto} B_1 \overset f_1 \to {\dotto} B_2 \dots
$$

This vertical set only depends on the underlying bidisks in a fairly crude way, as the following
Proposition shows.

\proclaim {Proposition \stabofcrosses}
If $U_m = U'_m \cup U''_m$ with $U'_m \cap U''_m \ne \emptyset$ and $U'_m$ and $U''_m$ homeomorphic
to disks. Set $B'_m=  U'_m \times V_m$ and $B''_m=  U''_m \times V_m$. Suppose
 $$
f_m: U_m \times V_m \to U_{m+1}\times V_{m+1}
$$
is an analytic map defined on an appropriate subset, such that the restrictions
$$
f'_m:B'_m \dotto B'_{m+1}\quad \text{and} \quad
f''_m:B''_m \dotto B''_{m+1}
$$
are crossed mappings of degree 1, then the stable sets of the sequences
$$
B'_0 \overset f'_0 \to {\dotto} B'_1 \overset f'_1 \to {\dotto} B'_2 \dots \quad \text{and} \quad
B''_0 \overset f''_0 \to {\dotto} B''_1 \overset f''_1 \to {\dotto} B''_2 \dots
$$
coincide.
\endproclaim

\demo{Proof}
In the proof of Corollary \unstablemansofcrosses\ above, the sequence $u_m$ could be chosen
arbitrarily, in particular in $U'_m \cap U''_m$.
\QED
\enddemo

\proclaim {Corollary \biinfinitesequences}
Let
$$
\dots B_{-1} = U_{-1} \times V_{-1}, B_0 = U_0 \times V_0, B_1 = U_1 \times V_1, \dots
$$
be a bi-infinite sequence of bidisks, and $f_i: B_i \dotto B_{i+1}$ be crossed mappings of degree 1,
with  $U'_i$ of uniformly bounded size in $U_i$. Then for all $m \in \Z$,
\newline
(1) the set
$$
\eqalign{
W^S_m = \{\,(x_m,y_m) \mid &\text{there exist $(x_n,y_n) \in B_n$}\cr
&\text{for all $n \ge m$ such that $f_n(x_n,y_n) = (x_{n+1},y_{n+1})$}\,\}}
$$
is a closed vertical-like Riemann surface in $B_m$, and $\pr_2: W^S_m \to V_m$ is an isomorphism;
\newline
(2) the set
$$
\eqalign{
W^U_m = \{\,(x_m,y_m) \mid &\text{there exist $(x_n,y_n) \in B_n$}\cr
&\text{for all $n < m$ such that $f_n(x_n,y_n) = (x_{n+1},y_{n+1})$}\,\}}
$$
is a closed horizontal-like Riemann surface in $B_m$, and $\pr_1: W^S_m \to U_m$ is an isomorphism.
\newline
(3) Moreover, the sequence
$$
(x_m,y_m) := W^S_m \cap W^U_m\,,\quad m\in\Z,
$$
is the unique bi-infinite sequence with $(x_m,y_m) \in B_m$ for all $m \in \Z$, and $f_m(x_m,y_m)
\allowbreak = (x_{m+1},y_{m+1})$.
\endproclaim

\demo{Proof}
The first statement is immediate from Corollary \unstablemansofcrosses, and the second also by
considering the mappings $g_n=f_{n+1}^{-1}$, which also define a bi-infinite sequence of crossed
mappings by Proposition \compsofcrosses (a). The third part follows immediately from the first two.
\QED
\enddemo

%% file: hen2_perturbpols.tex
\heading
\perturbpolschap. Perturbations of hyperbolic polynomials
\endheading

Let $p(z)$ be a hyperbolic polynomial of degree $k \ge 2$, which will be fixed for the next three
sections. We will drop the subscript $p$, and write
$$
H_a\left(\bvec xy \right) = \bvec{p(x) - ay}{x}.
$$ 

Choose as in section \telescopesection\ a neighborhood $U$ of $J_p$ such that $p:U' = p^{-1}(U) \to
U$ is a covering map. Set $U'' = p^{-1}(U')$.

Recall that when $a = 0$, the H\'enon mapping
$$
H_0: \bmatrix x\\y \endbmatrix \mapsto \bmatrix p(x) \\ x \endbmatrix
$$
maps all of $\C^2$ to the curve $C_p$ of equation $x = p(y)$, and reduces to $x \mapsto p(x)$ in the
first coordinate. Thus we can think of the polynomial $p$ as a mapping $C_p \to C_p$; when we think
of $U$ as a subset of $C_p$, we will denote it by $\tilde U$, and its projection onto the $y$-axis
simply by $U$.

First let us recall the crudest properties of H\'enon mappings; we will suppose $|a| \le 1$. If $p(z)
= a_kz^k + \dots + a_0 = a_kz^k + q(z)$, denote by $|q|(r) = |a_{k-1}|r^{k-1}+ \dots + |a_0|$, and
let $R$ be the largest root of the equation $|a_k|r^k - |q|(r) -2r =0$. We will call
$$
B_R = \altset {\bvec xy \in \C^2}{|x| < R, |y| < R}.
$$
All of the interesting dynamics of $H_a$ occurs in $B_R$, because Figure \crudepic\ roughly describes
the orbits of points.

%\setEPSF EPSF.crudepic 2.33in by 2.46in caption (Figure \crudepic: Crude
%picture of the dynamics of $H_a$.)
\setPSFig 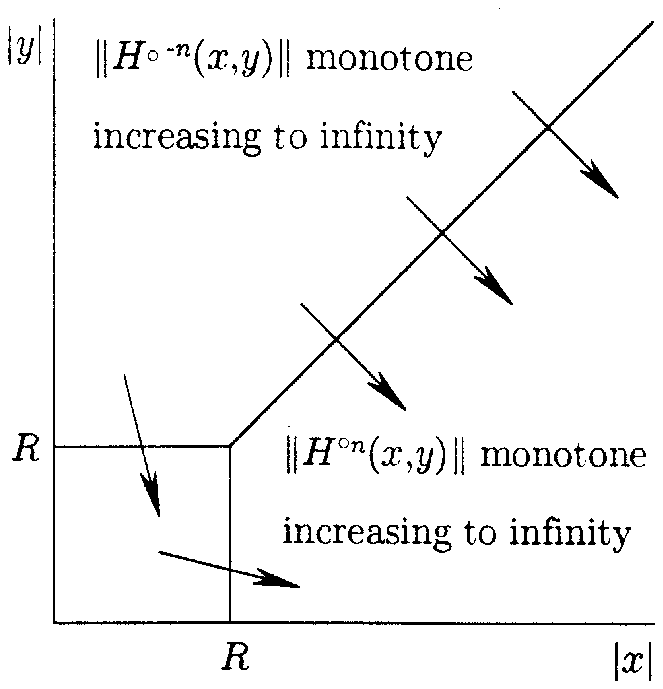 {3 true in} caption (Figure \crudepic: Crude
picture of the dynamics of $H_a$.)

Our construction will depend on two numbers $\delta > 0$ and $A > 0$, which will be chosen to satisfy
Requirements 1, 2, 3, 4, and 5, which are given below. Consider the neighborhood
$$
N_\delta = \altset{\bvec{x}{y} \in \C^2}{|p(y) - x| < \delta}
$$
of $C_p$.

% \setEPSF EPSF.deltaneib 2.24in by 1.54in caption (Figure \deltaneib: The
% curve $C_p$ and its neighborhood $N_\delta$, drawn in $\R^2$.)
\setPSFig 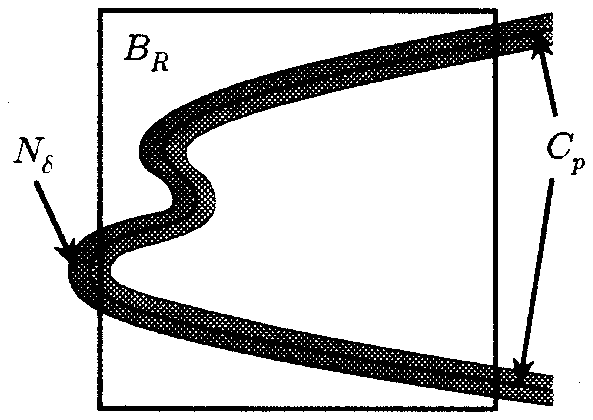 {2 true in} caption (Figure \deltaneib: The
curve $C_p$ and its neighborhood $N_\delta$, drawn in $\R^2$.)

Our first requirement concerns only $\delta$.

{\bf Requirement 1.}
The number $\delta>0$ is sufficiently small that $N_\delta$ intersects the boundary of $B_R$ only in
the ``vertical'' boundary $|x| = R$, and moreover for any $x_0 \in U$, each component of the
intersection $L_{x_0} \cap N_\delta$, where $L_{x_0}$ is the vertical line of equation $x = x_0$
contains a unique point of $C_p$, which will belong to $\tilde U'$. We will further require that
$$
\left|p\left(\zeta+\frac z{p'(\zeta)}\right)-(p(\zeta)+z)\right| < \frac \delta 2
$$
for all $\zeta \in J_p$ and $|z| < \delta$.

Choose, for the rest of the paper, a number $\delta$ satisfying Requirement 1.

{\bf Requirement 2.}
Now choose a number $\epsilon >0$ such that the sets $U_z, z \in U$ are all homeomorphic to disks, as
in Section \telescopesection. In Section \stablemansection, we will require a bit more: for all
$z \in J_p$, the image $p(U_z) \subset D_{\delta/2}(p(z))$, is contained in the Euclidean disk of
radius $\delta/2$ centered at $p(z)$. This will clearly be the case if $\epsilon$ is sufficiently
small.

Our next requirements all concern the size of $|a|$.

{\bf Requirement 3.}
We have $H_a(B_R) \subset N_\delta$ when $|a| < A$.

This will clearly be satisfied as soon as $A$ is sufficiently small.

Let 
$$
V' = \pr_1^{-1}(U) \cap N_\delta.
$$
be the union of these components. There is a well-defined function $u: V' \to U'$ given by $u(x,y) =
p^{-1}(x)$, the branch of the inverse image being precisely the intersection with $C_p$ above, which
one can also understand as the branch ``close to $y$''.

%\setEPSF EPSF.PictureOfV 4.28in by 4.00in caption (Figure \pictureofV: The
%neighborhood $V'$ of the Julia set $\tilde J_p \subset C_p$)
\setPSFig 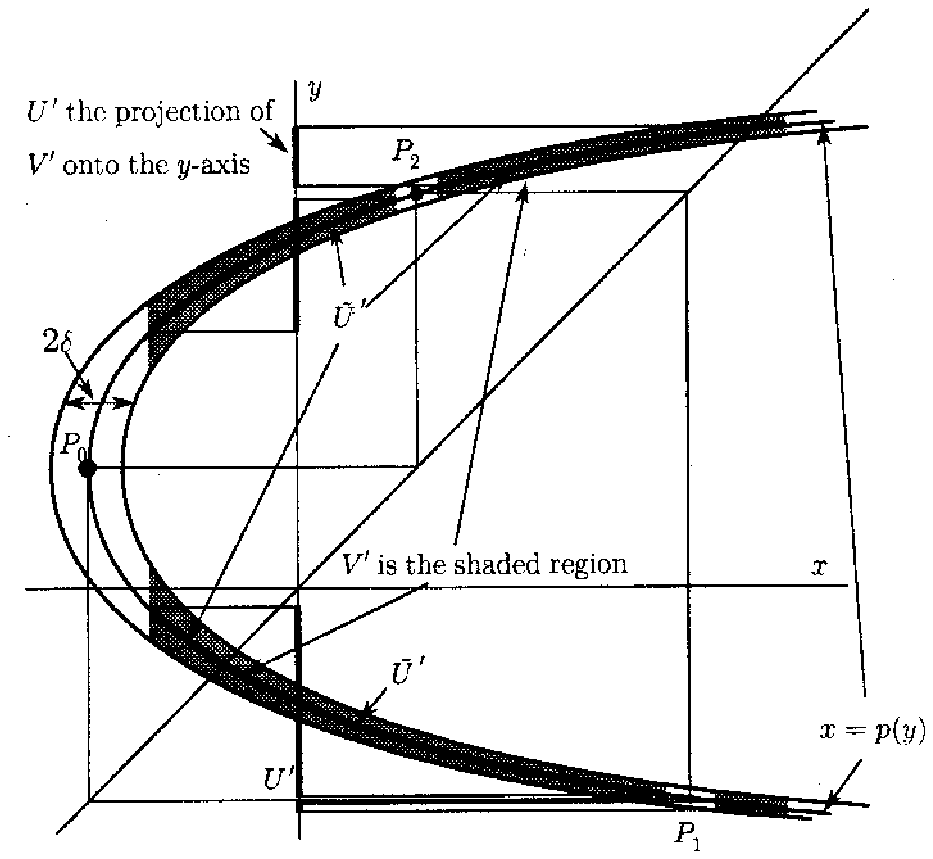 {4.75 true in} caption (Figure \pictureofV: The
neighborhood $V'$ of the Julia set $\tilde J_p \subset C_p$)

The pair of functions $(u,v): V \to \C^2$ given by the formulas
$$
u(x,y) = p^{-1}(x), \quad v(x,y) = p(y) - x
$$
parameterize $V$.

{\bf Requirement 4.}
We will require that $H_a$ should map the vertical boundary of $V'$ outside of $V$ when $|a|\le A$.

Again this will occur whenever $|a|$ is sufficiently small.

\proclaim{Proposition \attcyclespersist}
(a) For each attractive periodic point $z_0$ of $p$, there is an analytic function $z(a)$ defined for
$|a|<A$, such that $z(0)=z_0$ and $z(a)$ is an attractive cycle of $H_a$.
\newline
(b) The points of compact components of $\overline N_\delta \setminus V'$ are attracted to these
cycles, and the points of the unique non compact component iterate to infinity.
\endproclaim

\demo {Proof}
The union of the compact components of $\overline N_\delta \setminus V'$ is mapped into itself
by Requirement 3. Thus the sequence of iterates of $H_a$ is normal. On the other hand any limit
function has compact image. So the sequence of iterates is accumulating on finitely many attracting
cycles. The proof shows that these depend analytically on $a$ for $|a|<A$.
\QED
\enddemo

For every $z\in U'$ consider the neighborhood
$$
V_z = \altset {\bmatrix x\\y\endbmatrix\in V' }{u(x,y) \in U_z}
$$
of the point $(p(z),z) \in \tilde U'$.

\proclaim {Lemma \Visbidisk}
Under Requirement 3, the mapping
$$
\bvec xy \mapsto \bmatrix u(x,y\\ v(x,y)\endbmatrix
$$
is a biholomorphic isomorphism of $V_z$ onto the bidisk $U_z\times D_\delta$.
\endproclaim

The proof is left to the reader.

For all $z \in U''$, set
$$
W_z = V_z \cap H_a^{-1}(V_{p(z)}) \quad \text{and}\quad \tilde W_z = V_{p(z)} \cap H_a(V_z).
$$

\proclaim {Proposition \howintersect}
There exists $A>0$ such that if $|a|<A$, then for all $\zeta \in U''$, the mapping $H_a: W_z
\to\tilde W_z$ is a crossed mapping $V_z \dotto V_{p(z)}$ of degree 1.
\endproclaim

\demo {Proof}
Choose $\zeta\in \C$ with $|\zeta|<\delta$, and consider
$$
u \circ H_a: V_z \cap \{p(y) - x = \zeta\} \to \C.
$$
The disk $V_z \cap \{p(y) - x = \zeta\}$ is parameterized by $y$, and when the Jacobian $a$ of $H_a$
is zero, this map is simply $\zeta \mapsto p(y)$, and in particular maps the boundary of $U_z$
strictly outside $U_{p(z)}$, with degree $1$. This remains true for a sufficiently small
perturbation, in particular for $|a|<A$ when $A>0$ is small enough, and it is easy to see that if
$a$ is sufficiently small, then this will be true for all $z \in U''$ and $\zeta$ with $|\zeta|\le
\delta$. It follows that for such sufficiently small $A$, condition 1 of Definition \crossdef\
of a crossed mapping is satisfied, and the first half of condition 3.

For condition 2 and the second part of 3, we use the inverse mapping. For any fixed $z \in U''$
and $w \in p\left(U_{p(z)}\right)$, consider the vertical disk 
$$
\Delta_w=\left\{\bvec wy \in V_{p(z)}\right\};
$$
which the coordinate function $v$ maps isomorphicaly to the disk of radius $\delta$: on
its boundary, we have $p(y) - w = \delta$. Let us compute:
$$
v\left(H_a^{-1}\left(\bvec wy \right)\right) = p\left(\frac {p(y) - w}a\right)-x.
$$
This takes $\d \Delta_w$ to a large curve when $|a|$ is small, since $|(p(y) - w)/a| = \delta/|a|$
is large, and $p$ takes large values there.
\QED
\enddemo

\remark {Remark \degoneremark}
We do not need to calculate the degree of this mapping restricted to $\tilde W_{z}\cap
\Delta_w$. By Proposition \degsofcrosses, this must be one. It is not obvious from our
computation: we might rather have expected $k = \deg p$. This is because the set
$$
\altset {\bmatrix w\\y \endbmatrix \in \Delta_w}{|p(y) - w| < \delta}
$$
has $k$ components, one for each inverse image of $p(z)$ under $p$. Figure \discwithsmalldiscs\ %
illustrates this phenomenon.
\endremark

%\setEPSF EPSF.diskswithlit 4.65in by 2.01in caption (Figure
%\discwithsmalldiscs: How the image of a vertical disc intersects $N_\delta$)
\setPSFig 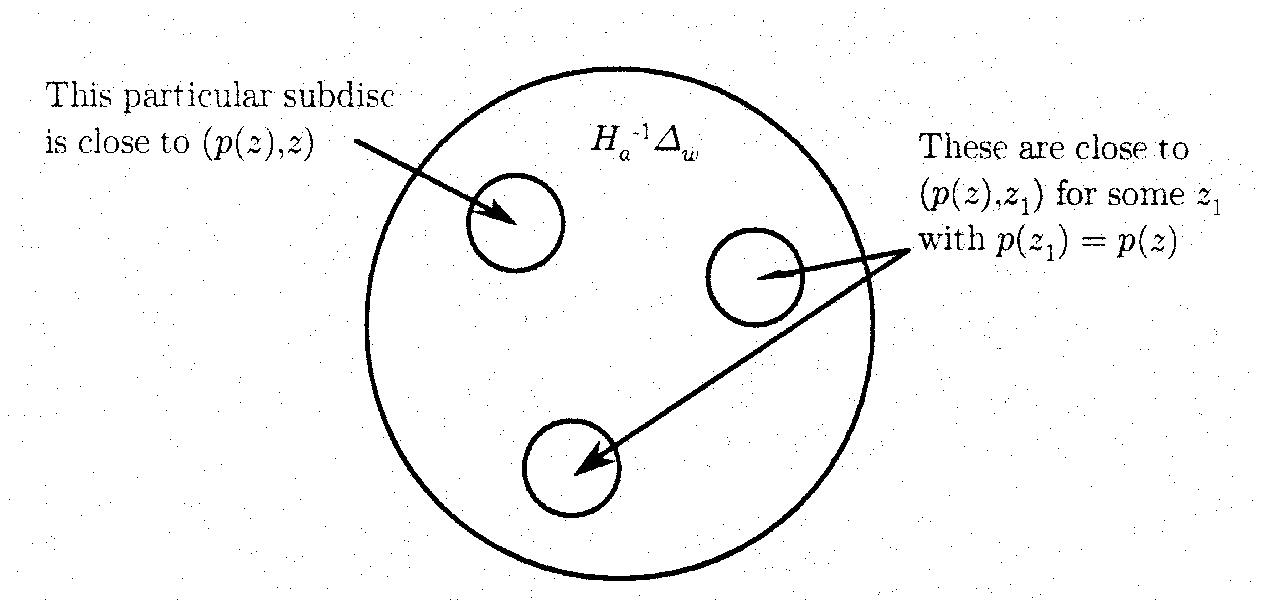 {2.5 true in} caption (Figure \discwithsmalldiscs: How
the image of a vertical disc intersects $N_\delta$) 

{\bf Requirement 5.}
The number $A$ is sufficiently small that the conclusion of Propo\-si\-tion \howintersect\ is
satisfied when $0<|a|<A$.

For any point $\underline z = (\dots, z_{-2},z_{-1},z_0)\in \hat J_p$, consider the bi-infinite
family of crossed mappings which we will denote, by abuse of notation
$$
\dots \overset H_a \to {\dotto} V_{z_{-1}} \overset H_a \to {\dotto} V_{z_0} \overset H_a \to
{\dotto} V_{p(z_0)} \overset H_a \to {\dotto} V_{p^{\circ 2}(z_0)} \overset H_a \to {\dotto} \dots.
$$
We can now define the mapping $\Phi$: by Proposition \onecrossescontract, there is a unique point
$\Phi(\underline z)\in V_{z_0}$ such that
$$
H_a^{\circ m}(\Phi(\underline z)) \in
\cases
V_{z_{-k}}           &\quad\text{for $m \le 0$}\\
V_{p^{\circ k}(z_0)} &\quad\text{for $m > 0$}
\endcases.
$$

\proclaim {Theorem \structureofJ}
The mapping $\Phi: \hat J_p \to \C^2$ is a homeomorphism onto $J$ which conjugates $\hat p$ to
$H_a$. \endproclaim

\demo {Proof}
The mapping $\Phi$ is obviously continuous. We will construct an inverse $\Psi: J \to \hat J_p$.

Observe first that $J \subset V$. Indeed, $J \subset B_R$, hence $J \subset N_\delta \cap B_R$. But
by Proposition \attcyclespersist, we know that the points in $N_\delta \cap B_R$ and not in $V$ tend
to $\infty$ or attracting cycles, hence cannot be points of $J$.

Therefore a point $(x,y) \in J$ defines a bi-infinite $p$-telescope $U_{u(H_a^{\circ m}(x,y))}$. We
have seen that such a bi-infinite telescope defines a point of $\underline z \in \hat J_p$. The
mapping
$\Psi: (x,y) \mapsto \underline z$ is obviously continuous. We leave it to the reader to check that
it is an inverse of $\Phi$.
\QED
\enddemo

This construction endows $J$ with a ``stable'' and ``unstable'' manifold: set $W^U$ (respectively
$W^S$) to be the union of all the unstable (respectively, stable) disks of the families of
1-crossed mappings defining $\Phi$.

\proclaim {Proposition \whereJs}
(a) We have the equalities
$$
J_+ \cap V = \bigcap_{n \ge 0} H_a^{-n}(V) \quad \text{and} \quad
J_- \cap V = \bigcap_{n \ge 0} H_a^{n}(V).
$$
\newline
(b) Moreover, the forward orbit of any point in $J_+$ is eventually contained in $V$, and the
backwards orbit of any point in $J_-$ is also eventually contained in $V$, except for the attracting
cycles described in Proposition \attcyclespersist.
\endproclaim

\demo{Proof}
The inclusions
$$
\bigcap_{n \ge 0} H_a^{-n}(V) \subset K_+ \cap V \quad \text{and} \quad
\bigcap_{n \ge 0} H_a^{n}(V) \subset K_- \cap V = J_- \cap V
$$
are obvious.

To see $\bigcap_{n \ge 0} H_a^{-n}(V) \subset J_+ \cap V$, we need to know that no interior point of
$K_+$ can have its forward orbit entirely in $V$. This follows from Proposition \stabofcrosses. On
such an open set, the sequence $H_a^{\circ n}, n \ge 0$ is normal, and we can extract a subsequence
$H_a^{\circ n_i}$ convergent on compact subsets. Given any two points $(x_1,y_1)$ and $(x_2,y_2)$,
and extracting a further subsequence if necessary, the infinite sequences of bidisks
$$
V_{u(H_a^{\circ n_i}(x_1,y_1))} \quad \text{and} \quad V_{u(H_a^{\circ n_i}(x_2,y_2))}
$$
connected by the mappings $H_a^{\circ (n_{i+1}-n_i)}$ are equivalent, hence define the same stable
set, which is a vertical-like disk in $V_{u(H_a^{\circ n_1}(x_1,y_1))}$. This contradicts the
assumption that
$(x_1, y_1)$ and $(x_2,y_2)$ could be chosen in an open set.

To see the opposite inclusions, consider a point of $V$. If its forward orbit ever leaves V, we know
what it does: it is either attracted to one of the attracting cycles described in Proposition
\attcyclespersist, in which case it is in the interior of $K_+$, or it iterates to infinity, in which
case it is not in $K_+$. Thus $J_+ \cap V = \bigcap_{n \ge 0} H_a^{-n}(V)$.

Now, suppose a point in $V$ has an inverse image not in $V$. Then since $H_a(N_\delta \cap B_R)
\subset N_\delta$, the inverse image is not in $N_\delta$, and a further inverse image is not in
$B_R$. By the discussion of Figure \crudepic, this means that the backwards orbit of the point tends
to infinity. This proves (a).

For (b), consider the forward image of any point in $K_+$. It will eventually be contained in
$N_\delta\cap B_R$, and if it is not eventually in $V$, then it is attracted to an attracting cycle,
by Proposition \attcyclespersist. Similarly, the backwards orbit of any point in $J_-$ will eventually
be contained in $N_\delta \cap B_R$. If it is in the basin of attraction of one of the attracting
cycles of Proposition \attcyclespersist, then its backwards orbit will enter and remain in $V$ unless
it is one of the attracting cycles itself.
\QED
\enddemo

\demo{Remark}
Proposition \whereJs\ is where we eliminate the possibility of wandering domains. In one dimension,
we have Sullivan's {\it No Wandering Domains}\/ Theorem to eliminate this possibility, but this uses
quasi-conformal mappings, and there does not appear to be an analog in several dimensions. For
hyperbolic polynomials (or rational functions), one can also eliminate the existence of wandering
domains by using the expanding metric on a neighborhood of the Julia set. The proof above is a
natural extension of that proof.
\enddemo

%% file: hen2_unstabman.tex
\heading
\unstablemansection. Characterization of $J_-$
\endheading

In this section we will prove the following result.

\proclaim {Theorem \unstabmanthm}
There is a homeomorphism
$$
\Phi_-: \hat\C_p \to J_-
$$
which conjugates $\hat p$ to $H_a$. This homeomorphism coincides with $\Phi: \hat J_p \to J$ on $\hat
J_p$.
\endproclaim

\demo{Remark \improveremark}
This result can be slightly improved: $\Phi_-$ can be chosen analytic on $\hat\C_p \setminus \hat
K_p$. But it cannot be made analytic on $\hat K_p$.
\enddemo

\demo{Proof}
We will begin by finding the restriction
$$
\Phi_-\mid_{\hat U'}: \hat U' \to W^U
$$
which conjugates $\hat p$ to $H_a$ there, then we will extend it to the remainder of $\hat \C_p$. This
requires the following.

\proclaim{Proposition \projtoU}
There exists a mapping $\pi_{U'}: V' \to U'$ which semi-conjugates
$$
H_a: H_a^{-1}(V') \cap V' = V'' \to V' \quad \text{to} \quad p: U'' \to U'.
$$
This mapping can be chosen so that for every $(x,y) \in V'$, we have
$$
\pi_{U'}(x,y) \in U_{u(x,y)},
$$
and so that the fibers are vertical-like.
\endproclaim

\demo{Proof}
We will begin by constructing our mapping on $V' \setminus V''$.

\proclaim{Lemma \existProjNearBound}
There exists a continuous mapping $\pi_{U'}: \overline{V'} \setminus V'' \to \overline{U'}
\setminus U''$ which semi-conjugates $H_a$ to $p$ as maps from the inner boundary to the outer
boundary, and such that the fibers are vertical-like disks. 
\endproclaim

We have found this lemma surprisingly difficult to prove. Note that $v: V' \setminus V'' \to
D_\delta$ is a locally trivial fibration, and our lemma says that there exists a trivialization with
special properties. Of course trivializations exist, since the base is contractible. But the
requirement that the induced sections be vertical-like does not seem to be accessible by
topological techniques. Instead, we will use differential equations. 

Before doing this, we require a fundamental statement about complex vector spaces.

\proclaim {Sublemma \hyperplanes} A real hyperplane $F$ in a complex vector space $E$ contains a
unique complex hyperplane $[F] = F \cap iF$.
\endproclaim

The statement contains the proof.

\demo{Proof of Lemma \existProjNearBound}
Consider the radial vector field $\xi = - r \partial
/\partial r$ on the disk $D_\delta$. This vector field lifts canonically to a vector field
$\Tilde{\Tilde \xi}$ on the ``vertical boundary'' of $\d^{\text{ver}}(\bar V' \setminus V'')$, both
inner and outer. This boundary is a 3-dimensional real manifold in $\C^2$, so at each point 
$$
(x,y) \in \d^{\text{ver}}(\bar V' \setminus V''),
$$
the tangent space 
$$
T_{(x,y)}\d^{\text{ver}}(\bar V' \setminus V'')
$$
contains a unique complex line
$$
[T_{(x,y)}\d^{\text{ver}}(\bar V' \setminus V'')],
$$
which maps isomorphically to $\C$ by the derivative $d_{(x,y)}v$. The vector field $\Tilde{\Tilde
\xi}$ is the unique lift of $\xi$ to this bundle of complex lines.

Let $\tilde \xi$ be a $C^\infty$ lifting of $\xi$ to $\bar V' \setminus V''$, which extends
$\Tilde{\Tilde \xi}$, and everywhere points into the vertical cone. Such a lifting exists, since
local liftings exist, and can be patched together by partitions of unity.

Denote by $\phi_{(x,y)}(t)$ the solution curve of the differential equation defined by $\tilde
\xi$ with $\phi_{(x,y)}(0)=(x,y)$. 

\proclaim {Sublemma \troublewithdiffeq}
The lift $\tilde \xi$ can be chosen so that
\newline
(a) The limit $w(x,y) := \lim_{t \to \infty} \phi_{(x,y)}(t)$ of this solution exists.
\newline
(b) The set $w^{-1}(\tilde z)$ is a vertical-like ``disk'' for all $\tilde z = (p(z),z) \in \tilde
U' \setminus \tilde U''$, although this disk {\it is not\/} an analytic disk in general unless $z$
is in the boundary $\d (\bar U'-U'')$.
\endproclaim

\demo {Proof of Sublemma \troublewithdiffeq}
We will work in the real oriented blow-up of $V'$ along $\tilde U'$. This is a set in which every
point of $\tilde U'$ is replaced by a circle. It is easy to describe in this case. Consider $\tilde
V' = U' \times ([0,\delta) \times S^1)$; the map $\tilde V' \to V'$ which sends $(z, r,\theta)$ to the
point $(x,y) \in V'$ with $u(x,y) = z$ and $v(x,y) = re^{i\theta}$ realizes $\tilde V'$ as such a
blow-up. In $\tilde V'$, we can lift not just $-r \partial/\partial r$ but also $-\partial/\partial r$,
and the existence and uniqueness theorem applies even to points on the boundary $U' \times \{0\}
\times S^1$. Now our disks are precisely the solutions which end on a circle $\{z\} \times \{0\}
\times S^1$.
\QEDL{(Sublemma~\troublewithdiffeq)}
\enddemo

The disks which foliate the vertical boundary of $\bar V' \setminus V''$ are collapsed to points under
$w$, in fact the point at which such a disk intersects $C_p$.

Thus it is now enough to choose a homeomorphism of
$$
\overline{\tilde U'} \setminus \tilde U''\quad \text{with} \quad \bar U' \setminus U''
$$
which extends the identity on the outer boundary, and maps a point $(x,y)$ of the inner boundary to
the point in $p^{-1}(u(H_a(x,y)))$ which is close to $(x,y)$.
\QEDL{(Lemma~\existProjNearBound)}
\enddemo

To prove Proposition \projtoU, we need to extend $\pi_{U'}$, and there is an obvious way to do so on
$V' \setminus W^S$: define
$$
\pi_{U'}(x,y) = p^{\circ -N}(\pi_{U'}(H_a^{\circ N}(x,y)))
$$
where $N$ is defined so that $H_a^{\circ N}(x,y) \in V' \setminus V''$, and the branch of
$p^{\circ -N}$ being chosen recursively so that $p^{\circ -m}(\pi_{U'}(H_a^{\circ N}(x,y)))$ is the
inverse image of $p^{\circ -m+1}(\pi_{U'}(H_a^{\circ N}(x,y)))$ in $U_{p^{\circ N-m}(u(x,y))}$.

This defines $\pi_{U'}$ on $V' \setminus W^S$. On $W^S$ it is defined by the telescope construction.
We still need to show that $\pi_{U'}$ is continuous; this is only an issue at points of $W^S$.

If $(x_k,y_k)$ is a sequence in $V'$ approaching $(x,y) \in W^S$, then the number $N_k$ of moves it
takes to escape $V'$ tends to $\infty$. Then the point $\pi_{U'}(x_k,y_k)$ is in the intersection of
the nested sequence
$$
\bigcap_{m=0}^{N_k} p^{-m}U_{u(H_a^{\circ m}(x_k,y_k))}
$$
which is a set with diameter tending to $0$ as $k \to \infty$. Moreover, the sets
$$
\bigcap_{m=0}^{N_k} p^{-m}U_{u(H_a^{\circ m}(x_k,y_k))} \quad \text {and}\quad
\bigcap_{m=0}^{N_k} p^{-m}U_{u(H_a^{\circ m}(x,y))}
$$
are close for $k$ large (in the sense that the Hausdorff distance of their closures are close, for
instance). Moreover, the second intersection converges $\pi_{U'}(x,y)$.
\QEDL{(Proposition~\projtoU)}
\enddemo

We can now construct our mapping
$$
\Phi_-: \hat U \to W^U.
$$
Given a point $(\dots, z_{-1}, z_0) \in \hat U$, let us consider the intersection
$$
\bigcap H_a^m(V_{z_{-k}}),
$$
which we have seen in Corollary \biinfinitesequences\ is a Riemann surface isomorphic to a
horizontal-like disk, in fact a section of the projection $v: V_{z_0} \to U_{z_0}$. This section
will intersect $\pi_{U'}^{-1}(z_0)$ in a single point: this point is $\Phi_-(\dots, z_{-1}, z_0)$.

With this definition, $\Phi_-$ is obviously continuous. And it is easy to construct an inverse:
simply associate to $(x,y) \in W^U$ the sequence
$$
z_{-n}(x,y) = \pi_{U'}(H_a^{\circ -n}(x,y)),\quad n \in \N.
$$
This is clearly a point of $\hat U$, and $\Phi_-(\underline z)(x,y) = (x,y)$.

Notice that the construction of $\Phi_-$ did not use the hypothesis that $H_a$ was injective, and
works even if $a=0$. In fact, it is especially easy in that case, giving the canonical projection
$\hat U' \to U'$, perhaps perturbed by a homeomorphism of $U'$ which commutes with $p$. Only the
construction of the inverse, \ie, the statement that $\Phi_-$ is injective, really uses the fact
that $H_a$ is an automorphism.

Now to extend $\Phi_-$ to the remainder of $\hat \C_p$. Any sequence of inverse images under $p$
tends to the Julia set, except for those finitely many periodic sequences consisting entirely of
points of the attracting cycles. Take a sequence $(\dots, z_{-1}, z_0)$. If it is not one of these
exceptional periodic histories, there exists $N$ such that $(\dots, z_{-N-1},z_{-N}) \in \hat U'$.
Map such a point to
$$
\Phi_-(\dots, z_{-1}, z_0) = H_a^{\circ N}\Phi_-(\dots, z_{-N-1},z_{-N}).
$$
If the sequence $(\dots, z_{-1}, z_0)$ is periodic, consisting of points of an attracting cycle, map
it to the attracting periodic point of $H_a$ corresponding to $z_0$.

We need to check that this is continuous at these exceptional points. Suppose that $(\dots, z_{-1},
z_0)$ is such an exceptional point, and that $(\dots, z'_{-1}, z'_0)$ is close to it. Let $N$ be
the first index such that $z'_{-N} \in U'$; the number $N$ is large. The point $\Phi_-(\dots,
z_{-N-1},z_{-N})$ is in $V \setminus V''$, and we have seen that such points have the same fate
under $H_a$ as $z_{-N}$ has under $p$, whether tending to infinity or to an attracting cycle. In
fact, we have seen further that the distance of $H_a^{\circ N}(\Phi_-(\dots, z_{-N-1},z_{-N}))$ to
the attracting cycle is arbitrarily small, by a quantity which depends only on $N$. This proves
continuity.

The only thing remaining is to check that the mapping $\Phi_-: \hat C_p \to J_-$ is surjective. We
invite the reader to check that if $H_a^{\circ -n}(x,y)$ is bounded as $n \to \infty$, then either
the sequence is eventually contained in $V'$, or it is the orbit of a repelling cycle for
$H_a^{-1}$.
\QEDL{(Theorem~\unstabmanthm)}
\enddemo

%% file: hen2_stabman.tex
\heading
\stablemansection. Characterization of $J_+$
\endheading

We are going to show that for $|a|$ sufficiently small, the set $J_+$ is modeled on $\check \C_p$.

First, we need to know a little more about this space: $\check \C_p$ is foliated by surfaces, in fact
by Riemann surfaces isomorphic to $\C$, each of which is dense. For each $\zeta \in J_p$, let
$L_\zeta$ be the inductive limit of
$$
\{\zeta\}\times D \overset f_p \to{\hookrightarrow} 
\{p(\zeta)\}\times D \overset f_p \to{\hookrightarrow}
\{p^{\circ 2}(\zeta)\}\times D \overset f_p \to{\hookrightarrow}\dots
$$ 
which is an increasing union of discs.

\proclaim {Proposition \foliationofCCech}
Each $L_\zeta$ is a Riemann surface isomorphic to $\C$, and is dense in $\check \C_p$.  The foliation
is compatible with the dynamics in the sense that
$$
\check p(L_\zeta) = L_{p(\zeta)}.
$$
\endproclaim

\demo{Proof}
The space $L_\zeta$ is a Riemann surface since the inclusions
$$
\{\zeta\}\times D \overset f_p \to{\hookrightarrow} 
\{p(\zeta)\}\times D \overset f_p \to{\hookrightarrow}
\{p^{\circ 2}(\zeta)\}\times D \overset f_p \to{\hookrightarrow}\dots
$$ 
are analytic inclusions; it is simply connected since it is a union of discs. The union is isomorphic
to $\C$ by Proposition 7.3 of \cite{HO}. The formula $\check p(L_\zeta) = L_{p(\zeta)}$ is obvious,
and it implies the density as follows. 

If $p^k(\zeta_1) = p^k(\zeta_2)=\zeta$, then $L_{\zeta_1} = L_{\zeta_2}$ since both are
equal to $\check p^{-k}(L_\zeta)$. Thus given any point $\zeta$, we see that $L_{\zeta'} =
L_\zeta$ for all $\zeta' \in p^{-k}(p^{\circ k}(\zeta)), k = 1, 2,3,\dots$. But such points
$\zeta'$ are dense in $J_p$, so the leaf $L_\zeta$ is dense in $J_p \times D$, which it
intersects in infinitely many discs $\{\zeta'\} \times D$.  It is then easy to show that it
is till dense in all the $((J_p \times D),n)$ in the inductive limit describing $\check
\C_p$.
\QED
\enddemo

Now, let us see the third part of Theorem \maintheorem.

\proclaim {Theorem \stabmantheorem}
(a) If $p$ is a hyperbolic polynomial and $|a|$ is sufficiently small, there exists a homeomorphism
$$
\Phi_+: \check \C_p \to J_+
$$
conjugating $\check p$ to $H_a\mid_{J_+}$.
\newline
(b) The mapping $\Phi_+$ maps the leaves of the foliation of $\check \C_p$ to Riemann surfaces
isomorphic to $\C$ immersed in $\C^2$.
\endproclaim

\demo{Proof}
We will construct $\Phi_+$ first on $J_p \times D = (J_p \times D) \times 0 \subset \check
\C_p$, and even then we will define it on larger and larger parts. The restriction of $\Phi_+$ to
$J_p \times D$ and all its further restrictions will be called $F_+$.

The first step (and the hardest) is to get started. We have already constructed $J$ and its stable
manifold $W^S$.

\proclaim {Lemma \stabmantojulia}
There is a unique projection $\pi: W^S \to J_p$ such that the diagram
$$
\CD
W^S      @>H_a>>   W^S    \\ 
@V\pi VV       @V\pi VV \\ 
J_p      @>p>>   J_p
\endCD 
$$
commutes, and the fibers of $\pi$ are stable disks of the crossed mappings.
\endproclaim

\demo{Proof of Lemma \stabmantojulia}
This is precisely what telescopes are for. Given a point $(x,y) \in W^S$, we can consider the
$p$-telescope $U_{u(H_a^{\circ k}(x,y)}$, which defines a unique point $\pi(x,y) \in J_p$.
\QEDL{Lemma~\stabmantojulia}
\enddemo

\proclaim {Proposition \strucOfWs}
There exists a homeomorphism
$$
F_+: J_p \times D \to W^S
$$
such that the following diagram commutes. 
%\setCD stabmanCD 2.00in by 1.50in

\centerline{\psfig{figure=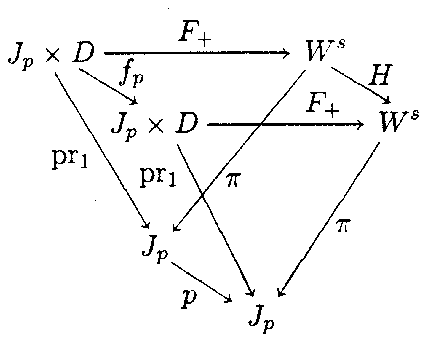}}
\endproclaim

This is a fairly difficult result and will require several steps, of which the first is the hardest.

\proclaim {Proposition \homeoonpartescaping}
There exists a homeomorphism 
$$
F_+: J_p \times \bar D \setminus f_p(J_p \times D) \to \bar W^S \setminus H_a(W^S)
$$
such that the diagram above commutes on the boundary.
\endproclaim

\demo {Proof of Proposition \homeoonpartescaping}
Before embarking on the proof, we will want to write our mapping $f_p$ in a different way, more
convenient for the rather delicate argument in Lemma \asymdev. Temporarily  write 
$$
g_p(\zeta,z) = \left(p(\zeta), a\left(\zeta+\frac z{p'(\zeta)}\right)\right), \quad\text{where}\quad
|z|<\delta.
$$
This mapping is conjugate to our original one: set $w = Rz/\delta$, so that in the coordinate $w$
we find
$$
(\zeta,w) \mapsto \left(p(\zeta), \frac {aR}\delta \left(\zeta + \frac {\delta w}{R p'(\zeta)}\right)
\right).
$$   
Thus if we choose $R = \delta/a$, we find that this linear change of variables conjugates the two
mappings.

Note that both $J_p\times D$ and $W^S$ are trivial bundles of disks of radius
$\delta$ over $J_p$. They are canonically homeomorphic, although that homeomorphism is not
compatible with the dynamics and is not the one we are trying to find. But we will use it
to define  
$$
F_+: J_p\times \d D \to \d W^S
$$
to be the ``identity'' (\ie, the map which takes $(\zeta,z)$ to the unique point $(x,y)$
above $\zeta$ such that $p(y) - x = z$). We can immediately define $F_+$ on the interior
boundary 
$$
F_+: f_p(J_p \times D) \to H_a(\d W^S)
$$
by the formula
$$
F_+ = H_a \circ F_+\circ f_p^{-1}.
$$
To extend it to the region in between, we will use some heavy-duty topology.  We define the space 
$$
\text{Homeo}_{J_p}(J_p \times D \setminus \check p(J_p \times D), W^S \setminus H_a(W^S);F_+)
$$
to be the fiber bundle over $J_p$, the fiber of which over $\zeta \in J_p$ is the space of
homeomorphisms of the fiber of $J_p \times D \setminus \check p(J_p\times D)$ above $\zeta$ to
the fiber of  $W^S \setminus H_a(W^S)$ above $\zeta$, which agree with $F_+$ (as defined so far)
on the boundaries.

\remark {Remark \infdimconst}
This space is pretty wild: it is a fiber bundle over $J_p$, the fiber of which is an
infinite-dimensional space of homeomorphisms of a surface with boundary, in fact a disk with holes,
to another. More precisely, the fiber of $J_p \times D \setminus f_p(J_p\times D) \to J_p$ is the
disk of radius $\delta$ with $k$ disks of radius $a\delta$ removed:
$$
D \setminus \bigcup_{\zeta_1 \in p^{-1}(\zeta)} D_{a\delta}(a\zeta_1).
$$
The fiber of $W^S \setminus H_a(W^S) \to J_p$ is similar.
\endremark

%\setEPSF EPSF.homeopic 4.86in by 3.07in caption (Figure \homeopic: The two fiber bundles)
\setPSFig 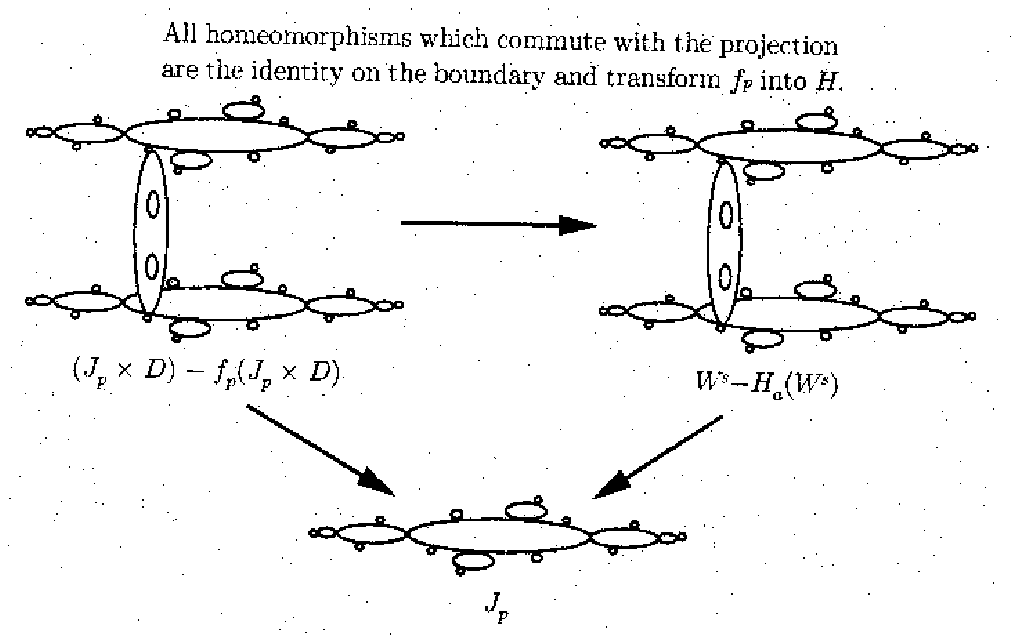 3.2in caption (Figure \homeopic: The two fiber bundles)

The statement of Proposition \homeoonpartescaping\ says exactly that the fiber bundle
$$
\text{Homeo}_{J_p}(J_p \times D \setminus f_p(J_p \times D), W^S \setminus H_a(W^S);F_+)
$$
has a continuous section. As a first step to proving this, we will require the following Lemma. 

\proclaim {Lemma \Hamstrom}
The fibers of the fiber bundle
$$
\roman{Homeo}_{J_p}(J_p \times D \setminus f_p(J_p \times D), W^S \setminus H_a(W^S);F_+)
$$
are contractible.
\endproclaim

\demo{Proof of Lemma \Hamstrom}
This is an immediate consequence of a hard theorem from topology, due to Hamstrom \cite{Ha}. This
theorem asserts that if $S$ is a compact surface with non-empty boundary, then the components of the
group of homeomorphism which are the identity on the boundary are contractible.  In our case, we can
choose a homeomorphism 
$$
D \setminus \bigcup_{\zeta_1 \in p^{-1}(\zeta)} D_{a\delta}(a\zeta_1) \to \pi_{U'}(\zeta) - H_a(W^S)
$$
which extends $F_+$ on the boundary by the classification of surfaces. Then the space of
homeomorphism homotopic to that one is acted on freely by the group of homeomorphisms of the
domain which are the identity on the boundary.
\QEDL {Lemma~\Hamstrom}
\enddemo

\remark{Remark \quasiconformal}
We are making this assertion for homeomorphisms. It is also true of $C^\infty$ diffeomorphisms \cite
{EE}, and there might be good reasons to prefer them. We are going to pull back the mapping $F_+$
repeatedly by diffeomorphisms, and eventually extend to a Cantor set. The resulting mapping will not
be differentiable on the Cantor set, but it will be quasi-conformal. The mapping $F_+$ would then be
a quasi-conformal map from the leaves of the foliation of $\check \C_p$ to the leaves of $J_+$,
showing that these Riemann surfaces are isomorphic to $\C$.  We will get this result by other means;
but it might be nice to know that as dynamical systems $\check \C_p$ and $J_+$ are quasi-conformally
isomorphic, in the sense of Riemann Surface Laminations.
\endremark

It follows that it is enough to show that there is section of the covering space of connected
components of the fibers, or alternately, that it is enough to show that there is a preferred
homotopy class of homeomorphisms from fibers to fibers, and there is: those which are homotopic
to a homeomorphism close to the identity. This requires a bit of elaboration.

\proclaim {Lemma \asymdev}
In the coordinates above, $H_a$ is close to $f_p$ when $a$ is small.
\endproclaim

\demo{Proof of Lemma \asymdev}
Let us start with a point $(\zeta,z) \in J_p \times D$, and consider $H_a^{-1} \circ f_p(\zeta,z)$,
where implicit in the composition is the parameterization of $W^S$ by the $(u,v)$-coordinates. Let
$(x_1,y_1)$ be the point of $W^S$ with $u(x_1,y_1) = p(\zeta)$ and $v(x_1,y_1) =
a(\zeta+z/p'(\zeta))$.

This has the very pleasant consequence that
$$
H_a^{-1}\left(\bvec {x_1}{y_1}\right) = \bvec {x_2}{y_2}
$$
where the $y2$-coordinate can be expressed exactly:
$$
y_2 = \frac1a (p(y_1) - x_1) = \dfrac1a v(x_1,y_1) = \zeta + \frac z{p'(\zeta)}.
$$
Thus 
$$
v(x_2,y_2)= p\left(\zeta +\frac z{p'(\zeta)}\right) - x_2 \approx z
$$
and we need to evaluate how good the approximation is. There are two approximations to consider:
$p(\zeta)= x_2 + \text{err}_1$ and $p\left(\zeta +\frac z{p'(\zeta)}\right) = p(\zeta)+z
+\text{err}_2$. 

We have that $|\text{err}_2| < \delta/2$ by Requirement 1 of Section \perturbpolschap. Moreover,
we have  $u(x_2,y_2) \in U_\zeta$, so $x_2 \in U_{p(\zeta)}$, and by Requirement 2 of Section
\perturbpolschap, this is smaller than $\delta/2$ 
\QEDL{Lemma~\asymdev}
\enddemo

Thus we can consider those homeomorphisms 
$$
h: J_p \times D\setminus f_p(J_p \times D) \to W^S \setminus H_a(W^S)
$$
which coincide with $F_+$ on the boundary, and which are homotopic with boundaries fixed to
homeomorphisms $h'$ which move all points by at most $|a|\delta$.  All such homeomorphisms are
homotopic, and this set is not empty, so it does define a homotopy class of homeomorphisms as
required.
\QEDL{Proposition~\homeoonpartescaping}
\enddemo

The next step is much easier.

\proclaim {Proposition \easyextension}
The mapping $F_+$ extends to a homeomorphism
$$
F_+: J_p \times D \setminus \hat J_p \to W^S \setminus J.
$$
\endproclaim

\demo {Proof of Proposition \easyextension}
Just map any point on $J_p \times D \setminus \hat J_p $ backwards by $f_p$ until it is in the
region $J_p \times D - f_p(J_p\times D)$, then map over by $F_+$ and back inwards by iterating
the H\'enon $H_a$ mapping the same number of times.
\QEDL{Proposition~\easyextension}
\enddemo

\proclaim {Proposition \phipluscont}
The map $F_+$ defined above and the map $\Phi$ from Theorem \structureofJ\ together give a
homeomorphism 
$$
F_+: J_p \times D \to W^S.
$$
\endproclaim

\demo {Proof of Proposition \phipluscont}
The only thing to prove is continuity, which is an issue only at points of $\hat J_p$. 
Suppose a sequence $(\zeta_n,z_n)$ in $J_p \times D$ converges to a point $(\zeta_0,z_0) \in \hat
J_p$.  Then the point $\zeta_0$ comes with a history, namely the first coordinates of 
$$
f_p^{\circ-1}(\zeta_0, z_0), f_p^{\circ-2}(\zeta_0, z_0),\dots
$$
which determines the point $(\zeta_0,z_0)$.  Clearly the points $F_+(\zeta_n,z_n)$ have
long backwards orbits $H_a^{\circ -m}(F_+(\zeta_n,z_n))$ contained in $W^S$, and the histories
$$
\pi_{U'}(H_a^{\circ -m}(F_+(\zeta_n,z_n))
$$
are close to that of $\zeta_0$.

The result follows from Proposition \onecrossescontract, which says that these long backwards
orbits restrict $z_n$ to a small disk, and Proposition \stabofcrosses, which says that this disc
is close to $z_0$.
\QEDL{Proposition~\phipluscont}
\enddemo

We can now complete the proof Theorem \stabmantheorem. Any point $((\zeta,z),n) \in \check \C_p$ has
a forward image in $J_p \times D$, namely 
$$
((\zeta,z),0) = \check p^{\circ n}((\zeta,z),n).
$$
So simply define 
$$
\Phi_+((\zeta,z),n) = H_a^{\circ -n}\left(F_+(\zeta,z)\right).
$$
This clearly defines an injective mapping $\Phi_+:\check \C_p \to J_+$.  It is surjective because
$J_+ \cap N_\delta = W^S$ by Proposition \whereJs.
\QED{Theorem~\stabmantheorem}
\enddemo

Finally, we will prove the following result, promised in the introduction.

\proclaim {Proposition \CChechWellDefined}
For all sufficiently small $\alpha_1$ and $\alpha_2$ and for all sufficiently large $R_1$ and $R_2$,
there is a homeomorphism 
$$
\psi: J_p \times D_{R_1} \to J_p \times D_{R_2}
$$
conjugating $f_{p,\alpha_1,R_1}$ to $f_{p,\alpha_2,R_2}$.
\endproclaim

\demo{Proof}
This is a consequence of Theorem \stabmantheorem: We proved that $H_a: W^S \to W^S$ is conjugate to
$f_{p,\alpha,R}: J_p\times D_R \to J_p \times D_R$  independent of $\alpha$ and $R$, so long as they
are respectively sufficiently small and sufficiently large. This certainly shows that 
$f_{p,\alpha_1,R_1}:J_p\times D_{R_1} \to J_p \times D_{R_1}$  and $f_{p,\alpha_2,R_2}:J_p\times
D_{R_2} \to J_p \times D_{R_2}$ are conjugate to each other under the same requirements.

It is also possible to prove this directly, by a proof very analogous to that of Theorem
\stabmantheorem. But we do not think that there is a much easier proof. Indeed, each of 
$f_{p,\alpha_i,R_i}:J_p\times D_{R_i} \to J_p \times D_{R_i}$ has infinitely many cycles,
each of which has multipliers in the vertical direction which depend on $\alpha_i$.  Thus
any conjugating map cannot be differentiable, so cannot be ``given by a formula'', and
some construction involving infinite processes must occur.
\QED 
\enddemo

%% file: hen2_examples.tex
\heading
\examplesection. Examples
\endheading

\subheading{Examples of $\hat \C_p$}

It is rather hard to find any particular category that $\hat{\C}_p$ belongs to in general. However,
when $p$ is hyperbolic, it is not too difficult to understand its structure. Except at the periodic
histories corresponding to the attracting cycles, the canonical projection $\hat \C_p \to \C$ is a
``ramified covering lamination'', with fibers homeomorphic to Cantor sets, and ramified above the
forward orbits of the critical points. The exceptional points have neighborhoods homeomorphic to
cones over solenoids.

The following statement is a first step towards seeing this.

\proclaim {Proposition \projnocrit}
If $\Omega \subset \C$ is a subset which does not intersect the orbit of the critical point, then
$\pi_\Omega: \hat \Omega \to \Omega$ is a locally trivial fibration with fiber a Cantor set.
\endproclaim

\demo {Proof}
Let $\hat z = (\dots, z_{-2}, z_{-1}, z_0) \in \hat \Omega$ and choose a connected, simply connected
neighborhood $\Omega'$ of $z_0$ in $\Omega$. Then for each $i$ there exists a branch $g_i$ of
$p^{\circ -i}$ defined on $\Omega'$ with $g_i(z_0) = z_{-i}$. Then the mapping
$$
(\hat z, v ) \mapsto (\dots, g_2(v), g_1(v), g_0(v))
$$
is a homeomorphism $\pi_\Omega^{-1}(z_0) \times \Omega' \to \hat \Omega'$.
\QED
\enddemo

Recall from the introduction that a Riemann surface lamination is a Hausdorff space locally
isomorphic to a product of a topological space with a Riemann surface. So long as a critical point is
not periodic, $\hat \C_p$ is still a Riemann surface lamination above the orbit. The only problem is
that the projection to $\C$ is ramified there. Indeed, let $z$ be any point not on an attracting
cycle, and choose $k$ such that all the points $y \in p^{ -k} (z)$ contain no post-critical points.
If $p$ is hyperbolic, then such a $k$ exists unless $z$ is a point of an attracting cycle, since the
orbits of the critical points are all attracted to the attracting cycles. Choose a
connected, simply-connected neighborhood $\Omega$ of $z$ such that there are no post-critical points
in $\Omega\setminus\{z\}$.

Now consider the following commutative diagram.
$$
\CD
\hat\Omega @>\hat p^{\circ-k}>> \bigcup\limits_{\text{components of $p^{-k}(\Omega)$}}
\hat\Omega'\\
@V\pi_{\Omega}VV                            @VV\bigcup \pi_{\Omega'}V \\
\Omega     @<<p^{\circ k}<      \bigcup\limits_{\text{components of } p^{-k}(\Omega)}\Omega'
\endCD
$$
The top mapping is an isomorphism in the category of Riemann surface laminations, so at any history
of a point except an attracting periodic point, the set $\hat \C_p$ is a Riemann surface lamination.
And since
$$
\pi_{\Omega} = p^{\circ k} \circ \bigcup \pi_{\Omega'} \circ \hat p^{\circ-k},
$$
we see that $\pi_{\Omega}$ is a ``ramified covering mapping'' of Riemann surface laminations,
ramified no worse than $p^{\circ k}$.

The space $\hat \C_p$ is {\it not}\/ a Riemann surface lamination at the periodic histories of
attracting periodic points.

First, let us examine the case $p_0: z \mapsto z^n$. In that case, $\hat \C_{p_0}$ is the cone over
the solenoid $\Sigma_n = \projlim(S^1, z \mapsto z^n)$. Indeed, in polar coordinates $p_0$
decouples, and the history of the radius contains no more information than the radius, since every
positive number has a unique positive $n$th root. In particular, $\hat \C_{p_0}$ is not a Riemann
surface lamination near the cone-point $\underline 0 =(\dots, 0, 0)$.

Next, let $p$ be any polynomial with an attractive cycle; we will see that $\hat \C_p$ is not much
nastier than the case above. Since $\hat \C_p$ and $\hat \C_{p^{\circ k}}$ are canonically isomorphic
for any $k \ge 1$, we may assume that $z$ is an attracting fixed point of $p$. Let $\Omega$ be the
component of $\C\setminus J_p$ containing $z$; and denote by $n$ the degree of $p$ mapping $\Omega$
to itself.

Let
$\underline z = (\dots, z, z)$ be the fixed point of
$\hat p$ corresponding to
$z$.

The space $\hat \Omega$ is connected only if $p$ has an attractive fixed point which attracts all
the critical points of $p$. In that case $J_p$ is a Jordan curve bounding $\Omega$, as in the case
of $p_0$ above. In general, $\hat \Omega$ is not connected, and the components of $\hat
\Omega$ are labeled by the totally disconnected set (projective limit of finite sets)
$$
\projlim_m (\pi_0(p^{-m}(\Omega)), \pi_0(p)),
$$
where $\pi_0(X)$ (the zeroth homotopy set) is the set of connected components of $X$, and
$\pi_0(p):\pi_0(p^{-m}(\Omega)) \to (p^{-m+1}(\Omega))$ is the map induced by p. The component to
which $\underline w = (\dots, w_{-1}, w_{0}) \in \hat \Omega$ belongs is indexed by
$$
(\dots, [w_{-1}], [w_0]),
$$
where the $[w_{-i}]$ denotes the component of $p^{-i}(\Omega)$ containing $w_{-i}$. There is a
distinguished component
$$
\Hat{\Hat\Omega} = (\dots, [z], [z]).
$$
of $\hat \Omega$.

The same argument as above show that the other components are Riemann surface laminations, but
$\Hat{\Hat\Omega}$ is not; in fact it is a cone over a solenoid analogous to
the case of $p_0$. To state this precisely, let $\hat \Bbb D \subset \hat \C_{p_0}$ be the component
of $\hat \C_{p_0}\setminus \hat J_{p_0}$ consisting of the histories of points in the unit disk
($\Bbb D$ is the unit disk).

\proclaim {Proposition \BadAtAttCycle}
There is a homeomorphism
$$
\psi: \Hat{\Hat\Omega}\to \widehat {\Bbb D}
$$
conjugating $\hat p$ to $\hat p_0$.
\endproclaim

\remark {Remarks \BadExRemark}
(1) In particular, $\underline z$ cannot have a neighborhood which is a Riemann
surface lamination.

(2) We must have $\psi(\underline z) = \underline 0$ since these are the unique fixed points
of $\hat p$ in $\Hat {\Hat \Omega}$ and of $\hat p_0$ in $\hat D$.
\endremark

\demo{Proof of Proposition \BadAtAttCycle}
We will first construct the homeomorphism $\psi$ near the boundary of $\Omega$. Take the circle
$\Gamma$ of radius $R$ centered at $z$ in $\Omega$, for the Poincar\'e metric of $\Omega$, with $R$
sufficiently large that the disk bounded by $\Gamma$ contains all the critical values in $\Omega$.
Let $\Gamma' = p^{-1}(\Gamma) \cap \Omega$, and let $A_0$ be the closed annular region between
$\Gamma$ and $\Gamma'$. Similarly, let $B_0 = \set {w \in \C}{r \le |w| \le r^{1/n}}$ for some
$0<r<1$.

\proclaim {Lemma \HomeoBetweenAnnuli}
There exists a homeomorphism $\bar\psi_0: A_0 \to B_0$ conjugating $p$ to $p_0$ on the inner
boundaries. \endproclaim

The proof is left to the reader; the details are given in \cite {DH1}.

Next extend $\bar \psi_0$ to $A_1 = p^{-1}(A_0)$, $A_2 = p^{-1}(A_1)$, etc., until it is defined on
the part $A = \cup A_i$ of $\Omega$ outside of $\Gamma$. These extensions exist by the lifting
criterion for covering spaces. This extended $\bar \psi$ now defines a homeomorphism $\psi: \hat A_p
\to \hat B_{p_0}$.

Now take any point $\underline w \in \Hat{\Hat \Omega}$. If $\underline w \ne \underline z$, then
for some $m>0$ we will have $\hat p^{\circ -m}(\underline w) \in \hat A$, so we can define
$$
\psi(\underline w) = \hat p_0^m \circ \psi \circ \hat p^{\circ -m}(\underline w).
$$
Since the construction can be made reversing the roles of $p$ and $p_0$, the map $\psi$ above is a
homeomorphism.
\QEDL{(Proposition~\BadAtAttCycle)}
\enddemo

\subheading{Examples of $\check \C_p$}

We showed, in  Proposition 4.4 of \cite{HO}, that the inductive limit of
$$
f_{p_0}: S^1\times D \to S^1 \times D, \quad
f_{p_0}(\zeta,z) = (\zeta^d, \zeta-\alpha \frac z {\zeta^{d-1}})
$$
is a 3-sphere with a solenoid removed, and we identified the map $\check p_0: \check \C_{p_0} \to
\check \C_{p_0}$ as the restriction of a certain map $\tau_{d,0}$ from the 3-sphere to itself with
two invariant solenoids, one attracting and one repelling.

Moreover, we showed in Theorem 3.11 of \cite{HO} that the conjugacy class of any injective mapping
from the solid torus to itself, with appropriate contraction and expansion properties depended only
on its homotopy class. Thus we obtain the following result:

\proclaim {Proposition \JPlusisSphere}
If $p$ is a polynomial with an attractive fixed point which attracts all the critical points of $p$,
then $\check \C_p$ is a 3-sphere with a solenoid removed, and $\check p$ is conjugate to
$\tau_{d,0}$.
\endproclaim

This gives quite a complete understanding of $J_+$ for the small perturbations of such polynomials.
For other polynomials, even hyperbolic, the situation is more complicated. However, there still are
$3$-spheres contained in $\check \C_p$.

Let $p$ be a hyperbolic polynomial with an attracting cycle; as above, by considering an iterate of
$p$, we may assume that all attracting periodic points of $p$ are fixed. Let $z$ be such an
attractive fixed point, and let $\Omega$ be the component of $\C\setminus J_p$ containing $z$. Then
$p\mid_\Omega:\Omega \to \Omega$ is a proper mapping of some degree $k$.

The boundary $\d \Omega$ is a Jordan curve, which is mapped to itself with degree $k$, and there
exist (exactly $k-1$) homeomorphisms $\gamma_\Omega: S^1 \to \partial\Omega$ such that
$$
\gamma(z^k) = p(\gamma(z)).
$$
Call $f_{\gamma_\Omega}:S^1 \times D \to S^1 \times D$ the mapping
$$
(\zeta,z) \mapsto \left(\zeta^k,\gamma_\Omega(\zeta) - \alpha
\frac z{p'(\gamma_\Omega(\zeta))}\right)
$$
where $\alpha$ is chosen as in Lemma \pcheckOK.

Consider now the diagram
$$
\CD
J_p \times D    @>f_p>>   J_p \times D   @>f_p>>   J_p \times D     @>f_p>>\dots\\
@A\gamma_\Omega\times \id AA    @A\gamma_\Omega\times \id AA   @A\gamma_\Omega\times \id AA \\
S^1 \times D    @>f_{\gamma_\Omega}>>   S^1 \times D   @>f_{\gamma_\Omega}>>   S^1 \times
D     @>f_{\gamma_\Omega}>>\dots
\endCD
$$
Clearly this leads to an embedding
$$
\varinjlim (S^1 \times D, f_{\gamma_\Omega}) \hookrightarrow \check \C_p.
$$

Using Theorem 3.11 and Proposition 4.4 of \cite{HO}, we can easily see that
$$
\varinjlim (S^1 \times D, f_{\gamma_\Omega})
$$
is homeomorphic to the 3-sphere with a solenoid removed; call $(\d \Omega)\check{} $ its image in
$\check \C_p$. Just as $\d \Omega$ is the boundary of the immediate basin of attraction of an attracting fixed point for the polynomial $p$, we would like to
say that $\Phi_+((\d \Omega)\check{}\,)$ is the boundary of the basin of attraction of the
corresponding attracting fixed point of $H_a$. This is nonsense: the basin of attraction is dense in
$J_+$, and the boundary is all of $J_+$. This has been proved in full generality by Bedford and
Smillie \cite{BS2}, and independently by Fornaess and Sibony \cite{FS}, and can easily be
verified directly in our case.

But there is a refinement of the notion of boundary: the {\it accessible boundary}. Suppose $U$ is
an open subset of a topological space $X$. Then the accessible boundary of $U$ in $X$ is the set of
$x \in X \setminus U$ for which there exists a continuous mapping $\eta: [0,1] \to X$ such a $\eta(0)
= x$ and $\eta(0,1] \subset U$.

\example {Example \AccessBdyEx}
Let $X = [0,1]$ and $U = X \setminus C$, where $C$ is the standard Cantor set. Of course, $U$ is
dense, so the boundary of $U$ is all of $C$. But the accessible boundary is just the set of
endpoints of intervals of $U$: any mapping $\eta$ as above must have image in a single interval.
These two ``boundaries''are very different: for instance, the accessible boundary is countable and
$C$ is uncountable.
\endexample

\proclaim {Theorem \boundaryOfBasin}
If $|a|$ is small as in the Requirements of Section \perturbpolschap, then $H_{p,a}$ has an
attractive fixed point $z(a)$ corresponding to $z$, and the accessible boundary of its basin is
$\Phi_+((\d \Omega)\check{}\,)$.
\endproclaim

\demo {Proof}
In $N_\delta$ there is a projection map $\pi:W^S \to J_p$. Our analysis in Sections \perturbpolschap
\ and \stablemansection\ show that $\pi^{-1}(\d \Omega)$ is the boundary of the immediate domain of
$z(a)$ in $N_\delta$, and this agrees with the accessible boundary of that component. But
$\pi^{-1}(\d \Omega)$ is the first stage in the construction of $\Phi_+((\d \Omega)\check{}\,)$.

Now let $\eta:[0,1] \to K_+$ be an access to a point $(x,y) = \eta(0) \in J_+$ with $\eta(0,1]$ in
the basin of $z(a)$. Since the image of $\eta$ is compact, some forward image will be in
$\N_\delta$, and some further forward image will be in the component basin in $N_\delta$
containing $z(a)$, hence in $\Phi_+((\d \Omega)\check{}\,)$. But this set is stable under $H_{p,a}$,
so $(x,y) \in \Phi_+((\d \Omega)\check{}\,)$.
\QED
\enddemo

We will now describe what we know of the algebraic topology of $\check \C_p$. This is quite
difficult when $J_p$ is not a Jordan curve. For spaces like these, which are not locally
contractible, the only really well behaved theory is \v Cech cohomology; unfortunately, we will see
in Theorem \Jplussimplyconnected\ that $\check H^1(\check \C_p,G) = 0$ for all coefficients
$G$, and \v Cech cohomology carries no information. Thus we are forced to consider homology;
and there does not appear to be any really good homology theory. We will use singular
homology ``faute de mieux''; it may be that \v Cech or Steenrod homology are better behaved.

\example{Example \clamexample}
Consider the ``Hawaiian earring'' space
$$
X=\bigcup_{k\in \N} \altset {\bvec xy \in \R^2}{ x^2 + \left(y+\frac1k\right)^2 = \frac1{k^2}}.
$$
%\setEPSF EPSF.clam 2.01in by 2.19in caption (Figure \Hawaii: The Hawaiian earring)
\setPSFig 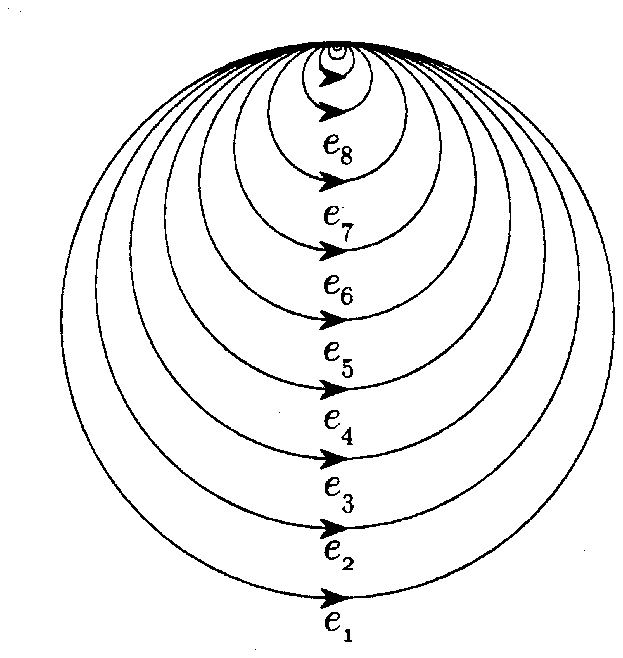 {2.75 true in} caption (Figure \Hawaii: The Hawaiian earring)

There is a canonical mapping $H_1(X,\Z) \to \Z^\N$ which associates to a cycle $\alpha$ the
sequence $\alpha(n), n \in \N$ where $\alpha(n)$ is the number of times $\alpha$ turns around the
$n$th circle. It is easy to show that the mapping is surjective, but it is apparently not
injective. At least, we do not see how a loop
$$
e_1e_2\dots e_1^{-1} e_2^{-1}\dots
$$
in the fundamental group can be written as a commutator.
\endexample

Suppose that $p$ is hyperbolic, that $J_p$ is connected, and that all the attractive cycles are
fixed. Call $\script X = \pi_0(\overset \circ \to {K_p})$, and $p_*: \script X \to \script X$ the map
induced by $p$. Let $\script X_0 \subset \script X$ be the finite subset of components containing
attractive fixed points. For each $X \in \script X$ the integer $k(X)$ is the degree
of $p$ restricted to the component $X$.

The space $J_p$ is very much like the Hawaiian earring, and there is an analogous mapping
$H_1(J_p,\Z) \to \Z^X$ which is surjective but probably not injective; nevertheless we consider the
kernel as pathological.

Recall from Corollary 4.5 of \cite{HO} that when $p$ is of degree $d$ and $J_p$ is a Jordan curve,
so that $\script X = \script X _0$ has a single element, then $\check \C_p$ is homeomorphic to the
complement of a solenoid in
$S^3$, and that
$$
H_1(\check \C_p,\Z) = \Z\left[\frac 1d\right].
$$
The same holds for the spaces $(\d X)\check{}$ for all $X \in \script X _0$, with $d$ replaced by
$k(X)$.

\proclaim {Theorem \homologyofbasin}
The inclusion
$$
\bigcup_{X \in \script X_0} (\d X)\check{} \to \check \C_p
$$
induces a split injection
$$
\bigoplus_{X \in \script X _0} \Z\left[\frac 1{k(X)}\right] \to H_1(\check\C_p,\Z).
$$
This mapping is not surjective.
\endproclaim

\demo{Proof} Consider the following diagram.
%\setCD 3dHomologyCD 4.18in by 1.72in

\centerline{\psfig{figure=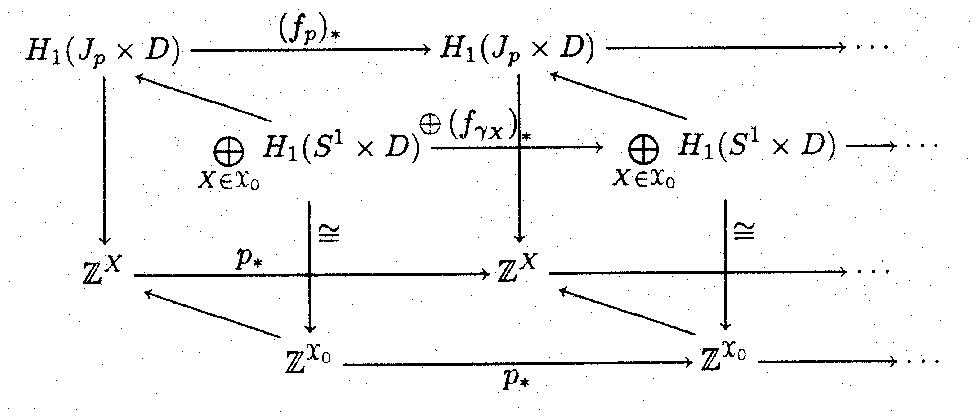}}

For each face
$$
\CD
\bigoplus_{X\in\script X_0} H_1(S_1 \times D) @>\bigoplus(f_{\gamma_X})_*>>  H_1(J_p \times D) \\
                            @VVV                                                    @VVV       \\
                      \Z^{\script X_0}                    @>>>                \Z^X
\endCD
$$
the horizontal mappings are split injections. This is obvious for the bottom mapping, and for the
top mapping the inclusion
$$
J_p \to J_p \cup \bigcup_{X \in \script X\setminus \script X_0} X
$$
induces a splitting on the homology. The theorem follows by passing to the direct limit, since the
direct limit is an exact functor, and the homology of the direct limit is the direct limit of the
homology.
\QED
\enddemo

Finally, why is \v Cech cohomology useless in this setting?

\proclaim {Theorem \Jplussimplyconnected}
If $J_p$ is not a Jordan curve, all covering spaces of $\check \C_p$ are trivial.
\endproclaim

\remark {Remark \BigFundGroupRemark}
It might seem that this is just another way of saying that $\check \C_p$ is simply connected, which
evidently contradicts Theorem \homologyofbasin. To resolve the apparent contradiction, notice that
$\check \C_p$ is not locally simply connected, so it doesn't have a universal covering space, and the
``singular'' fundamental group defined using loops does not classify covers. On the other hand, its
abelianization is the singular homology, so this fundamental group is enormous.
\endremark

\demo{Proof of Theorem \Jplussimplyconnected}
A covering $Z \to \check \C_p$ restricts to covers $Z_i \to (J_p \times D)\times \{i\}$, together with
inclusions $Z_i \subset Z_{i+1}$ such that the following diagram commutes.
$$
\CD
          Z_0             @>>>              Z_1              @>>>   \cdots \\
         @VVV                              @VVV                            \\
(J_p\times D)\times\{0\} @>>f_p>  (J_p\times D)\times\{1\}  @>>f_p> \cdots
\endCD
$$
Of course, the cover $Z_i$ of $J_p \times D$ restricts to a cover $Y_i$ of $J_p$, and the diagram
above gives covering homeomorphisms $\alpha_i: Y_i \to p^*Y_{i+1}$.
%\setCD covhomeoCD 2.65in by 0.67in

\centerline{\psfig{figure=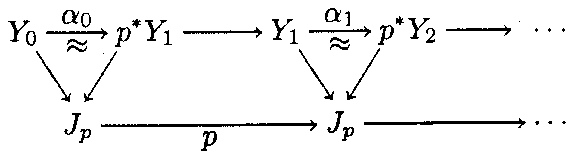}}

Now consider the cover $Y_0$: it is ramified at most over finitely many of the Jordan curves $\d
X, X \in \script X$. Indeed, you can find a cover $\Cal U = (U_j)$ of $J_p$ such that over every
$U_j$,$Y_0$ is trivial, and since $J_p$ is compact, we may take the cover finite. There then
exists a number $\delta$ for every $\zeta \in J_p$, the $\delta$-ball around $\zeta$ is
contained in one $U_j$. Since all but finitely many $X \in \script X $ have diameter smaller than
$\delta$, $Y_0$ is trivial over the boundaries of such components $X$.

Denote by $\script X '\subset \script X $ the set of components $X$ such that $Y_0$ is trivial over
$\d X$; notice that there exists $k_0$ such that $p^{\circ k}(X') = X$ for all $k \ge k_0$.

Choose $k \ge k_0$, any component $X \in \script X $, and $X' \in \script X '$ such that
$p^{\circ k}(X') = X$. If $Y$ is ramified above $\d X$, then $Y_0 \cong (p^{\circ k})^* Y_k$ is
ramified above $\d X'$. Since this is not the case, all $Y_k$ are unramified above all $\d X, X \in
\script X $. But this implies that $Y_k$ is trivial for all $k$, hence that all $Z_k$ are trivial,
hence that $Z$ is trivial.
\QEDL{(Theorem~\Jplussimplyconnected)}
\enddemo

\remark {Remark \SingvCechRemark}
This proof actually shows more than we claimed: it show that all principal $G$-bundles over $\check
J_p$ are trivial. Since the \v Cech cohomology $\check H^1(\check \C_p,G)$ classifies such principal
bundles, this shows that $\check H^1(\check \C_p,G) = 0$ for all coefficient groups $G$ (even
non-abelian, if you know how to define such things).
\endremark

It certainly seems remarkable that singular homology is picking up so much more than \v Cech
cohomology; one might expect the opposite. The following example should help to explain how this
happens, as well as give some insight into how $J_+$ is made.

\example {Example \hornexample}
Consider a tube $[0,\infty)\times S^1$ embedded in $\R^3$ so that it spirals onto a circle, as
suggested in Figure \horn.

%\setEPSF EPSF.horn 2.38in by 1.60in caption (Figure \horn: How a
%non-trivial circle may support no non-trivial covers)
\setPSFig 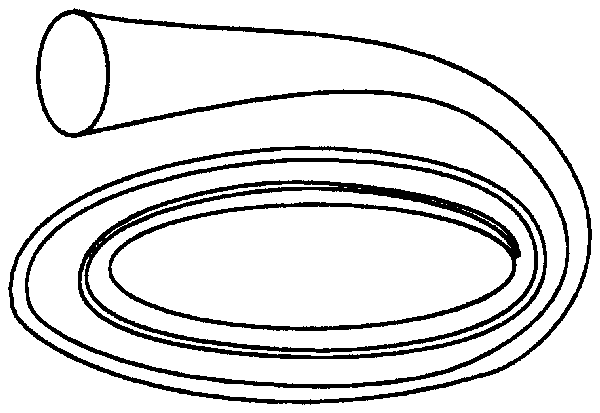 {2 true in} caption (Figure \horn: How a
non-trivial circle may support no non-trivial covers)

Let $C$ be the union of the circle and the tube. Then the singular homology is $H_1(C,\Z) =
\Z^2$, generated by the circle and the boundary of the tube. But no covering can be non-trivial
over the boundary of the tube, because it would then also be non-trivial over tiny cross-sections of
the tube near the circle, and such cross-sections will be contained within a single open set of any
open finite cover.
\endexample

This is the way $\check \C_p$ is made. There are big Jordan curves in $J_p$, but there are tubes in
$\check \C_p$ joining them to all their inverse images in $J_p$, which become arbitrarily small.

%% file: hen2_Wada.tex
\heading
\Wadasection. Lakes of Wada in Dynamical Systems
\endheading

A famous example in plane topology, due to Wada, is that there exist three bounded, connected and
simply connected open sets in $\R^2$ such that $\partial U_1 = \partial U_2 = \partial U_3$. We wish
to show that under appropriate circumstances the components of the basin of attraction of an
attractive cycle for a H\'enon mapping will form Lakes of Wada \cite{Y}.

The classical construction of Lakes of Wada illustrates the perils of philanthropy. Consider a
circular island, inhabited, to the sorrow of the others, by three philanthropists. One has a lake of
water, another of milk and a third of wine. The first, in a fit of generosity, decides to build a
network of canals bringing water within 100 meters of every spot of the island. It is clearly possible
to do this keeping the union of the original water lake and the water canals connected and simply
connected, with closures disjoint from the other lakes.

Next the second, perhaps worried about child nutrition, decides to bring milk to within 10 m of every
spot on the island, and builds a system of canals to that effect. She also keeps her milk locus
connected and simply connected.

Not to be outdone, the purveyor of wine now decides to bring wine to within 1 m of every spot on the
island. He finds his canal building rather more of an effort than the previous two, but being properly
fortified, he carries it out.

In turn, each of the three philanthropists brings his or her product closer to the poor inhabitants.
It should be clear that the construction can be continued, and that in the limit the construction
achieves the desired result: each of the lakes, being an increasing union of connected, simply
connected open sets, is a connected, simply connected set, and each point of the boundary of one is in
the boundary of the other two.

We will show that under appropriate circumstances, the basins of attraction of attracting cycles form
Lakes of Wada for H\'enon mappings in $\R^2$. As it turns out, the ``strategy'' of these basins is
remarkably similar to that of the philanthropists.

More specifically, we will work with {\it dense polynomials}. Let $p$ be a real hyperbolic polynomial
with connected Julia set, and suppose all the attracting cycles of $p^{\circ k}$ are real fixed
points. We will say that $p$ is dense if for each such fixed point $x$, its real domain of
attraction $\Omega_x \cap \R$ is dense in $J_p \cap \R$.

There are lots of dense polynomials. The following lemma describes some of them in degree 2. We have
found this lemma to be harder to prove than we had expected.

\proclaim {Lemma \goodpols}
Let $p$ be a real quadratic polynomial with an attracting cycle of period $k$, with $k$ an odd prime.
Then the $k$ basins $U_1=0, \dots ,U_{k-1}$ of the attracting fixed points of $p^{\circ k}$ in $\R$ are
all dense in $J_p \cap \R$.
\endproclaim

\demo{Proof}
Denote by $I_0$ the largest bounded interval invariant under the polynomial; it is bounded by the
``external'' fixed point and its inverse image. Without loss of generality we may assume that the
critical point is periodic of period $k$; let $c_0, c_1, \dots, c_{k-1}, c_k = c_0$ be the critical
orbit; all the interesting dynamics occurs in the interval $I = [c_1,c_2]\subset I_0$.

The polynomial $p$ also has an ``internal'' fixed point $\alpha \in [c_0,c_1]$. If $J \subset I$ is
any interval containing $\alpha$, then $\cup p^{\circ n}(J) = I$. The alternative is that $\cup
p^{\circ n}(J) = J_0$ is an interval in $[c_0,c_1]$ bounded by a cycle of period 2, and there are no
such cycles in $[c_0,c_1]$ (here we are using that $p$ is a polynomial, not just a unimodal map). It
follows from this that each of the basins $U_i$ accumulates at $\alpha$.

Thus to prove the lemma, it is enough to show that the real inverse images of $\alpha$ are dense in
the real Julia set $J_p \cap \R$. Let us denote by $V_0, \dots, V_{k-1}$ the immediate domains of
attraction in $\R$. It is known that if $k$ is an odd prime (or more generally simply odd) the $V_i$
have disjoint closures; let $\script T = \{T_1, \dots, T_{k-1}\}$ be the bounded components of $I
\setminus \cup V_i$. 

\proclaim {Sublemma \enoughtofindpointsinT}
If there is an inverse image of $\alpha$ in each $T_j$, then $p$ is a dense polynomial.
\endproclaim

\demo {Proof of \enoughtofindpointsinT}
The Julia set is 
$$
J_p \cap \R = I_0 \setminus \bigcup_{i=0}^{k-1} \bigcup_{n=0}^\infty p^{-n}(V_i).
$$
If each component of
$$
X_M =I_0 \setminus \bigcup_{i=0}^{k-1} \bigcup_{n=0}^M p^{-n}(V_i)
$$
contains an inverse image of $\alpha$, then these inverse images will accumulate on all of $J_p
\cap \R$. But if each component of $X_M$ contains an inverse image of $\alpha$, then this is also
true of each component of $X_{M+1}$, since $p$ maps each component of $X_{M+1}$ to a component of
$X_M$. Thus it is enough to start the induction, which is the hypothesis of the sublemma.
\QEDL{Sublemma~\enoughtofindpointsinT}
\enddemo

There is a repelling cycle $Z$ of length $k$ such that all endpoints of intervals $T \in \script T$
are either in $Z$ or in its inverse images. Let us denote $\script T'$ those intervals for which at
least one end-point is periodic, and $\script T''$ the others. Moreover set 
$$
A= \bigcup_{T \in \script T'} \bigcup _{n=0}^\infty p^n(T).
$$

Now there are two possibilities:

(a) If $\alpha \in A$, there is an inverse image of $\alpha$ in some $T'\in \script T'$. But then
there must be an inverse image of $\alpha$ in every $T \in \script T$, since each endpoint of $T$ will
eventually land on every point of $Z$, in particular on an end-point of $T'$; that iterate of $T$ will
cover $T'$. Then by sublemma \enoughtofindpointsinT, $p$ is dense.

(b) If $\alpha \notin A$, then A is disconnected, and $p$ permutes the components of $A$ circularly,
with period $k'$ with $1 < k' < k$. This is because some interval $T \in \script T'$ must have both
endpoints in $Z$, as there is one more point in $Z$ than there are intervals in $\script T$. That
interval must return to itself in fewer that $k$ moves. Moreover $k'$ divides $k$, since the map $Z
\to \pi_0(A)$ is equivariant, \ie, the following diagram commutes.
$$
\CD
  Z   @>\psi>>   \pi_0(A)    \\
@V p VV            @VV \pi_0(p) V  \\
  Z   @>>\psi>   \pi_0(A) 
\endCD
$$
This cannot happen if $k$ is prime.
\QEDL{(Lemma~\goodpols)}
\enddemo

Figure \periodninefig\ should illustrate what is going on.

%\setEPSF EPSF.periodnine 4.01in by 1.38in caption (Figure \periodninefig:
%The polynomial $z^2 - 1.785866 \dots$, with an attractive cycle of length
%9)
\setPSFig 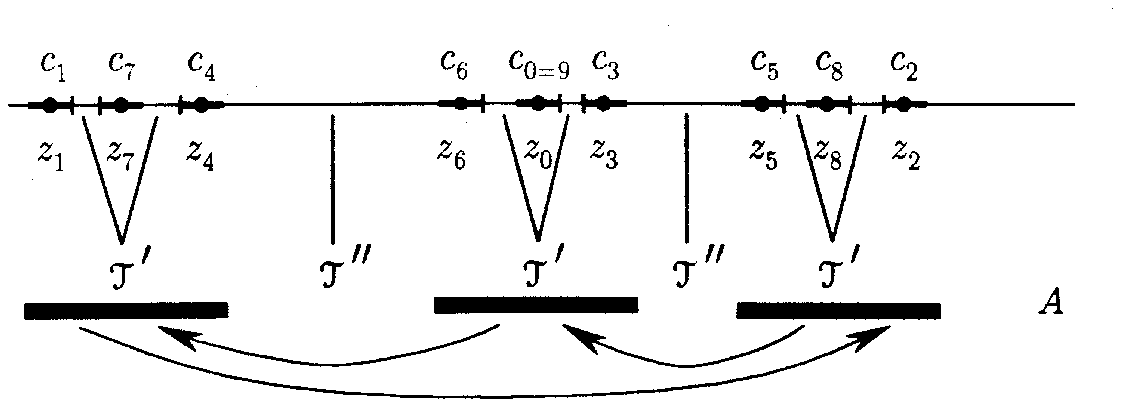 {1.75 true in} caption (Figure \periodninefig:
The polynomial $z^2 - 1.785866 \dots$, with an attractive cycle of length 9)

For this polynomial, the critical point is periodic of period 9. We have used heavy lines to
indicate the immediate basin, and the line segments pointing down form the repelling cycle $Z =
\{z_0, \dots, z_8\}$. The 8 intervals forming $\script T$ break up into 6 in $\script T'$, and two in
$\script T''$. The forward images of the intervals in $\script T'$ form the set $A$ which consists
of 3 intervals which are permuted circularly. The point $\alpha$ is not in $A$, and this polynomial
is not dense.

\remark{Remark \renormalizabilityremark}
The proof above should shows that if a hyperbolic polynomial is not dense, then it is renormalizable
in an appropriate sense. We could get necessary and sufficient conditions for a quadratic polynomial
to be dense by pushing the argument a bit further.
\endremark 

\proclaim{Theorem \wadalakestheorem}
If $p$ is a dense polynomial and if $|a|$ is sufficiently small, then the H\'enon mapping $H_{p,a}$
has attractive cycles close to those of $p$, and the boundaries of all the components of the basins
coincide.
\endproclaim

\remark {Remark \notcomplexgenerality}
General theorems of Bedford and Smillie \cite{BS3}, and independently by Sibony and Fornaess
\cite{FS}, assert that for any saddle point of a H\'enon mapping (and many other mappings besides),
the stable manifold is dense in $J_+$. We will use an analogous statement, in the much more
restricted class of mappings to which Theorem \stabmantojulia\ applies. But Theorem
\wadalakestheorem\ does not immediately follow from this density argument. For instance, the mapping
$$
\bvec xy \mapsto \bvec {x^2 - 1.05 - .38 y}x
$$
has an attractive cycle of period 3 (as well as an attractive fixed point), and the basin of this
cycle is bounded by the stable manifold of a cycle of period 3 which is a saddle. Of course, in
$\C^2$, each path component of this stable manifold is dense in $J_+$, and in particular each path
component accumulates onto the others. But not in $\R^2$: in the real, each of these path components
accumulates exactly on the stable manifold of the saddle fixed point.
\endremark

\demo{Proof}
The proof is contained in Lemma \stabmantojulia, Proposition \foliationofCCech\ and Proposition
\boundaryOfBasin. Let us review how these fit together to give the result.

Notice that the proof of Lemma \stabmantojulia\ is valid over the reals. Thus for $|a|$ sufficiently
small, $\Phi_+: \check \R_p \to J_+ \cap \R^2$ is a homeomorphism, where the space 
$$
\check \R_p= \varinjlim \left(J_p \cap \R)\times I , f_p\mid _{(J_p \cap \R)\times I}\right)
$$
is obtained by the same inductive limit construction as in the complex. Figures \LakesA, \LakesB,
\LakesC and \Lakes\ illustrate this construction.

Moreover, Proposition \boundaryOfBasin\ is also valid over the reals: if $x$ is a fixed point of
$p^{\circ k}$ with immediate basin $\Omega$, the accessible boundary of each basin is 
$$
(\d (\Omega \cap \R))\check{} = \varinjlim (\d \Omega\times I, f_p^{\circ k}).
$$
 
But $\Omega \cap \R$ is an interval, bounded by a repelling fixed point $\xi$ of $p^{\circ k}$ and one
of its inverse images $\xi'$. As such, the inductive limit above is a real line, which maps by
$\Phi_+$ to the the stable manifold of the fixed point $\xi(a)$ of $H_{p^{\circ k},a}$. Thus we
understand exactly what the accessible boundary of each basin is, and what its inverse image by
$\Phi_+$ is. So far, none of this required that $p$ be dense.

If $p$ is dense, then every point of $J_p \cap \R$ can be approximated by inverse images $\xi_n
\in p^{-nk}(\xi)$; the curves $\pi_{U'}^{-1}(\xi_n)$ are then part of $(\d (\Omega \cap \R))\check{}$,
by the argument of Proposition \foliationofCCech. Thus $(\d (\Omega \cap \R))\check{}$ is dense in
$(J_p \times I)\times \{0\}$, the first term in the inductive limit defining $\check \R_p$, and by the argument of \foliationofCCech, this shows it is dense in all of
$\check \R_p$. Thus the accessible boundary of each basin is dense in $J_+ \cap \R^2$, so they do have
common boundary.
\QED 
\enddemo

The following pictures carry out the construction of $\check \R_p$ for $p$ a real quadratic polynomial
with an attractive cycle of period 3. It is of course easy to imagine the first step of the
construction $(J_p \cap \R) \times I$, which is a product of a Cantor set by an interval.

%\setEPSF EPSF.check1 4.01in by 0.82in caption (Figure \LakesA: The set
%$(J_p\cap\R)\times I$; the first step in the construction)
\setPSFig 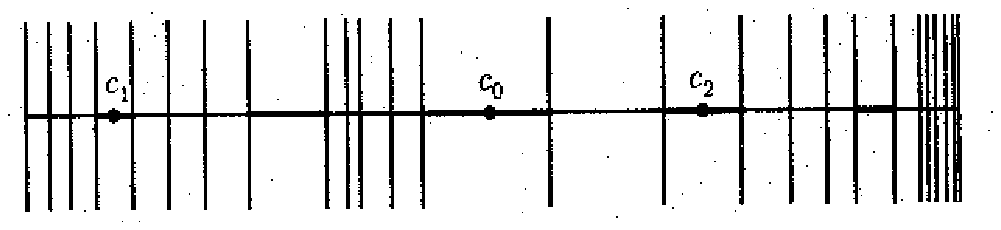 {1 true in} caption (Figure \LakesA: The set
$(J_p\cap\R)\times I$; the first step in the construction)

We have drawn a few genuine points of the Cantor set, and others ``impressionistically''. 

How should we imagine the inclusion 
$$
\left((J_p\cap\R)\times I\right) \times \{0\} \hookrightarrow
\left((J_p\cap\R)\times I\right) \times \{1\}?
$$
Note $f_p$ maps the two intervals through the endpoints of the immediate basin of $c_0$ to two
disjoint subintervals in the interval through the right endpoint of the immediate basin of $c_1$. Note
also that the $p'(\zeta)$ in the denominator in the definition of $f_p$ is essential for the
orientations to be as indicated by the arrows in Figure \LakesB.

%\settopEPSF EPSF.check2 4.01in by 2.13in caption (Figure \LakesB: The set
%$((J_p\cap\R)\times I)\times\{1\}$; the second step in the construction)
\setPSFig 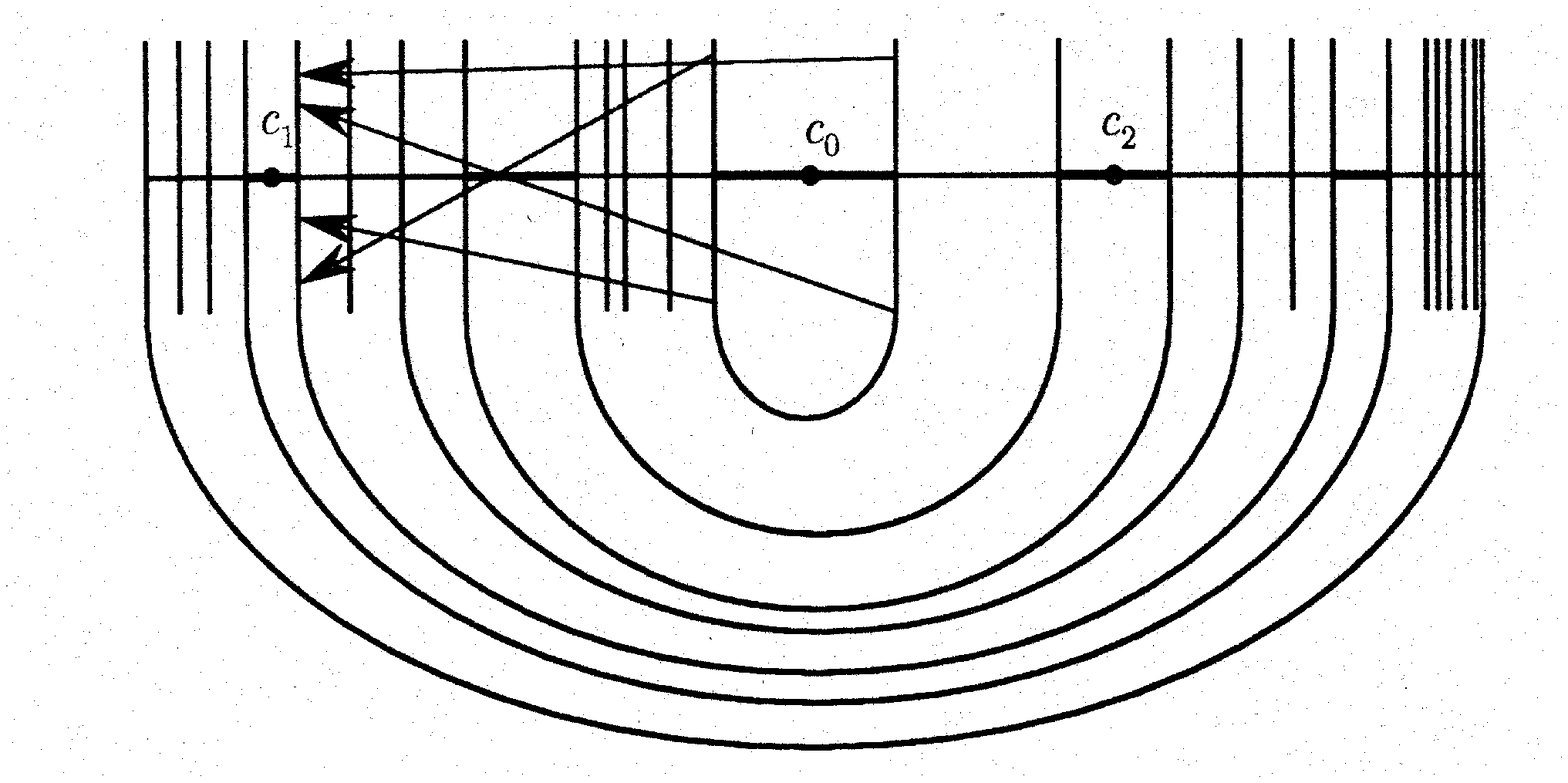 {6.5 true cm} caption (Figure \LakesB: The set
$((J_p\cap\R)\times I)\times\{1\}$; the second step in the construction)

Thus in $\left((J_p\cap\R)\times I\right) \times \{1\}$ there must be an arc joining the two intervals
above, so that these intervals and the arc will map to the interval where the arrows end. Similarly
one sees that there must be an arc joining every pair of symmetric intervals.

%\setEPSF EPSF.check3 4.01in by 2.86in caption (Figure \LakesC: The set
%$((J_p\cap\R)\times I)\times\{2\}$; the third step in the construction)
\setPSFig 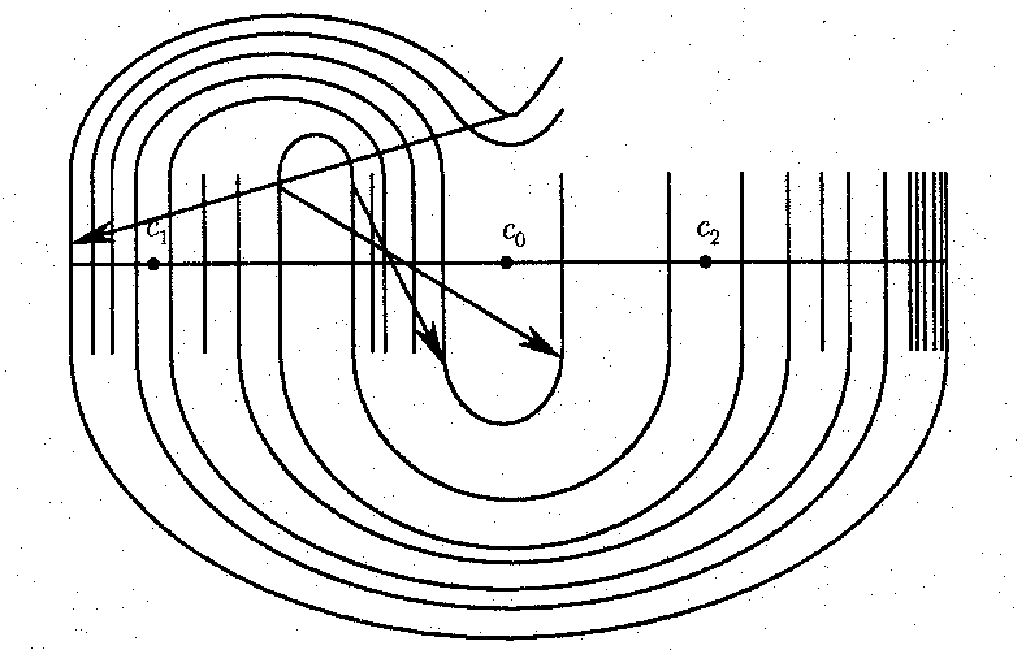 {3.5 true in} caption (Figure \LakesC: The set
$((J_p\cap\R)\times I)\times\{2\}$; the third step in the construction)

Figure \LakesB\ illustrates this construction. How should we continue the construction? In
$\left((J_p\cap\R)\times I\right) \times \{1\}$ we need inverse images of the arcs added in the
previous step; Figure \LakesC\ illustrates how this is to be done. Note that this time some of these
arcs do not join intervals to intervals. This is because points to the left of $c_1$ have no inverse
images in the Cantor set $J_p \cap \R$.

Making these pictures is a bit addictive, and if one gets carried away, the result may look like
Figure \Lakes.

%\settopEPSF EPSF.wadasmall 3.65in by 5.03in caption (Figure \Lakes: How the
% basins of the attractive cycle fit to form Lakes of Wada)
\setPSFig 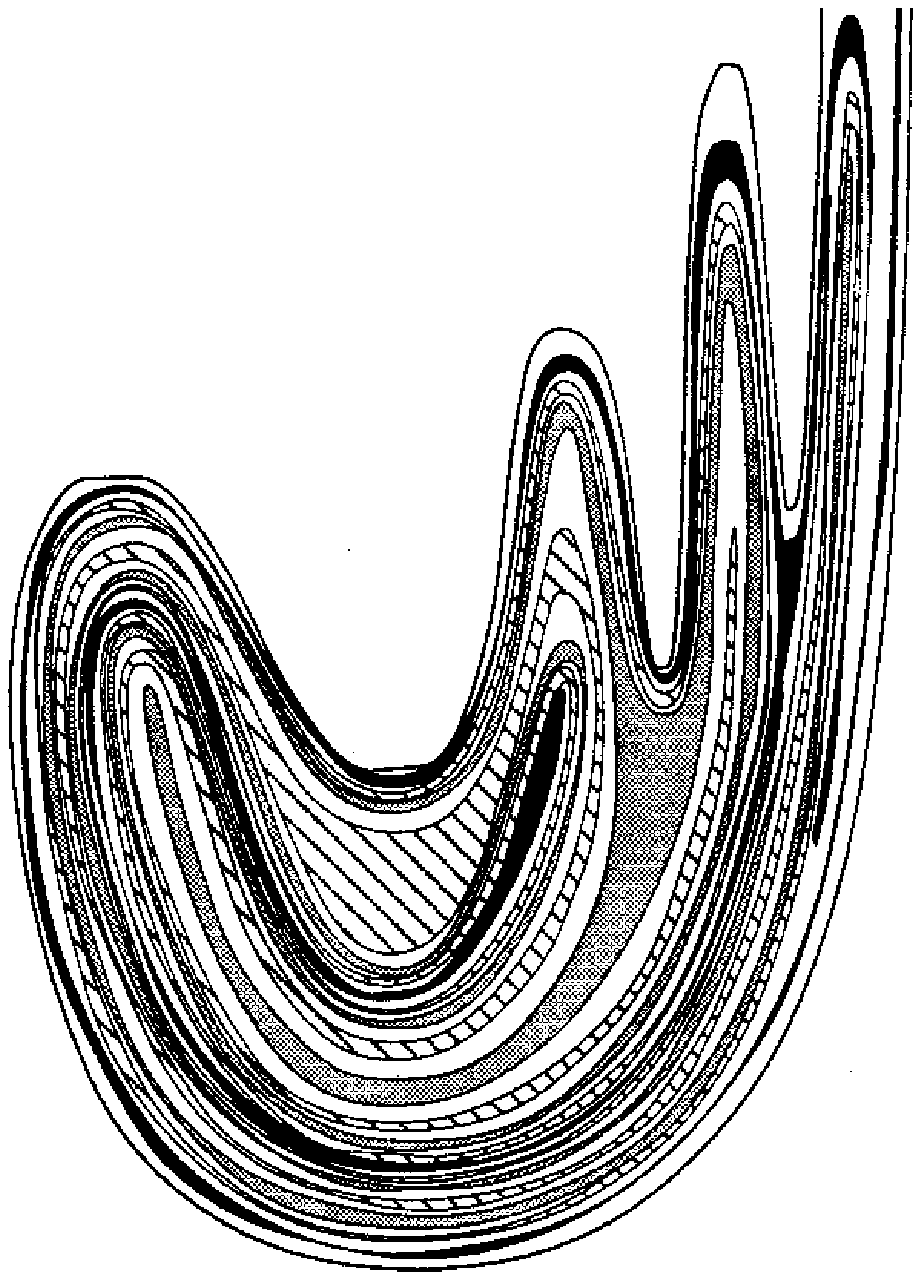 {6 true in}  caption (Figure \Lakes: How the
 basins of the attractive cycle fit to form Lakes of Wada)

\remark{Exercise for the Continuum Theorists}
Let $p$ be a dense quadratic polynomial, and denote by $X_{p,a}$ the one point compactification
of $J_+\cap \R^2$. Since this is naturally a subset of the sphere $S^2$, we can apply Alexander
duality, to get
$$
\check H^1(X_{p,a},\Z) = \overline H^0(S^2 \setminus X,\Z) \cong \Z^k
$$
if $p$ has an attracting cycle of period $k$. So for different $k$ these sets are certainly not
homeomorphic. What happens when $p_1$ and $p_2$ have the same period, but belong to different
components of the Mandelbrot set? The first interesting case is $k=5$, where there are 3 different
candidates. We would speculate that the space $X_{p,a}$ are not homeomorphic, but it isn't clear what
topological invariant distinguishes between them. Perhaps something could be made from the fact that
in the construction of the basins, the sequence of lefts and rights which the ``dead branches'' make
when they leave the ``main branch'' is different for these three polynomials, essentially reflecting
the kneading sequences.
\endremark
\bigskip

\subheading {Acknowledgments}
We wish to thank many people for helpful conversations, comments, and encouragement, including E.
Bedford, A. Douady, J.-E. Forn{\ae}ss, S. Friedland, D. Giarrusso, J. Mayer, C. McMullen, J. Milnor,
J. Rogers, N. Sibony, J. Smillie,  D. Sullivan, and J.-C. Yoccoz. We thank H. Smith for much computer
experimentation over the years.

We especially wish to thank Bodil Branner for organizing the conference in Hiller{\o}d, hosting a
magnificent party, and encouraging, cajoling, and helping us to finish this paper.